\magnification=1200

\centerline
{\bf Isometries of quadratic spaces}
\bigskip
\centerline
{Eva Bayer--Fluckiger}
\vskip 1,5cm

{\bf Introduction}
\bigskip
Let $k$ be a field of characteristic not 2. A {\it quadratic space} is a non--degenerate
symmetric bilinear form $q: V \times V \to k$ defined on a finite dimensional $k$--vector
space $V$, and an {\it isometry} of $(V,q)$ is an element of $O(q)$, in other words an
isomorphism $t : V \to V$ such that $q(tx,ty) = q(x,y)$ for all $x,y \in V$. In [M 69],
 Milnor raised the following question  : 
\bigskip
\noindent
{\bf Question 1.} Which quadratic spaces admit an isometry of given irreducible
minimal polynomial ?
\bigskip
The case of local fields is covered in [M 69], and the present paper gives an answer to
Milnor's question for global fields. 
\bigskip
Let $q$ be
a quadratic space, and let $f \in k[X]$ be an irreducible  polynomial.
The following Hasse principle is proved in section 9 (cf. th. 9.1) :
\bigskip
\noindent
{\bf Theorem.} {\it Suppose that  $k$ is a global field. The quadratic space $q$ has an isometry with
minimal polynomial $f$  if and only if such an isometry exists over  all
the completions of $k$.}
\bigskip
 In order
to obtain a necessary and sufficient criterion, we need to consider the case of reducible minimal polynomials
over local fields and the field of real numbers. This leads to a generalization of the above question.
Note that
any endomorphism $t : V \to V$ gives rise to a torsion $k[X]$--module. We ask the following :
\bigskip
\noindent 
{\bf Question 2.} Which quadratic spaces admit an isometry with a given torsion module ?
\bigskip

Note that this covers several special cases of interest :

\medskip
- The ``vertical" case : if $M  = [k[X]/(f)]^m$ where $f \in k[X]$ is an irreducible polynomial
and $m \in {\bf N}$
is precisely the question of Milnor mentioned above; 

\medskip
- The ``horizontal" case : if $M = k[X]/(f_1 \dots f_r)$ with $f_i \in k[X]$ distinct
irreducible polynomials, then question 2 amounts to asking which orthogonal
groups contain a maximal torus of a given type (see 
for instance [BCM 03], [G 05], [PR 10], [GR 12], [F 12], [Lee 12], [B 12]);

\medskip
- The case of ``rational knot modules", see for instance [Le 80].

\bigskip Integral analogs of this question arise in connection with algebraic-geometric and arithmetic applications (cf.
Gross-McMullen  [GMc 02], and  [BMa 13]).

\bigskip
Most of the  results concern fields of cohomological dimension 1, local and global fields; let
us  illustrate them by a few examples. Let $M$ be a self--dual torsion $k[X]$--module
with characteristic polynomial $F_M \in k[X]$
(see  \S 2), and suppose that ${\rm dim}(q) = {\rm dim}_k(M)$. We shall see that it
is sufficient to answer the question in the case of semi--simple modules
(cf. prop. 4.1). We
have (see cor. 6.3) :
\bigskip
\noindent
{\bf Proposition.} {\it Suppose that $k$ is a field of cohomological dimension $\le 1$, and that $M$ is semi--simple.
Then the quadratic space $q$ has an isometry with module $M$
if and only if ${\rm det}(q)F_M(1)F_M(-1) \in k^2$.}
\smallskip
\noindent
[Note that if $F_M(1)F_M(-1) = 0$, then this means that any quadratic space has an isometry with
module $M$, provided their dimensions are the same.]
\bigskip
In the case of global fields, we give an answer to Milnor's question.
Suppose that $f \in k[X]$ is 
an irreducible, monic polynomial such that $f(X) = X^{{\rm deg}(f)}f(X^{-1})$ and that $f \not = X  + 1$. Let $m \in {\bf N}$ and set $F = f^m$. Assume that ${\rm dim}(q) = {\rm deg}(F)$.
Then we have  (see cor. 9.2) :
\bigskip
\noindent
{\bf Theorem.} {\it Suppose that $k$ is a global field. The quadratic space $q$ has an isometry with minimal polynomial $f$
if and only if the signature condition and the hyperbolicity condition are satisfied {\rm (see
\S 9)}, and if ${\rm det}(q) = F(1)F(-1) \in k^*/k^{*2}$.}

\bigskip
The paper is structured as follows. The first three sections contain some definitions
and basic facts, including some results of Milnor [M 69]. Sections 4 and 5 are concerned with
isometries with a given module over an arbitrary ground field, and are used throughout
the paper. The following sections treat the case of fields of cohomological dimension 1
(\S 6), local fields (\S 7),  the field of real numbers (\S 8), and global fields (\S 9-\S12).

\bigskip
{\bf \S 1. Quadratic spaces, isometries and symmetric polynomials}
\bigskip
Let $k$ be a field of
characteristic not 2. 
A  {\it quadratic space} is a pair $(V,q)$, where $V$ is
a finite dimensional $k$--vector space, and $q : V \times V \to k$ is a
symmetric bilinear form of non--zero determinant. 
The {\it determinant} of  $(V,q)$ is
denoted by ${\rm det}(q)$; 
it is an element of $k^*/k^{*2}$. Let $n = {\rm dim}(V)$. The {\it discriminant}
of $q$ is by definition ${\rm disc}(q) = (-1)^{(n-1)n/2}{\rm det}(q)$. Any quadratic
space can be diagonalized, in other words there exist $a_1,\dots,a_n \in k^*$ such
that $q \simeq <a_1, \dots ,a_n>$. Let us denote by ${\rm Br}(k)$ the Brauer group of $k$,
considered as an additive abelian group, and let ${\rm Br_2}(k)$ be the subgroup of elements
of order $\le 2$ of ${\rm Br}(k)$. The {\it Hasse invariant} of $q$ is
by definition $w(q) = \Sigma_{i < j} (a_i,a_j) \in {\rm Br_2}(k)$, where $(a_i,a_j)$ is the class of the quaternion
algebra over $k$ determined by $a_i,a_j$. 
For more information concerning basic results on 
quadratic spaces, see for instance [O'M 73] and [Sch 85]. 
\bigskip
 An {\it isometry} of an quadratic space
 $(V,q)$ is an 
isomorphism $t : V \to V$ such that $q(tx,ty) = q(x,y)$ for all $x,y \in V$.

\bigskip
A monic polynomial $f \in k[X]$ is said to be
$\epsilon$--{\it symmetric}  for some $\epsilon = \pm 1$ if  $f(X) = \epsilon X^{{\rm deg}(f)} f(X^{-1})$. We say that $f$ is {\it symmetric} if it is $\epsilon$--symmetric for $\epsilon = 1$. 
The following is well--known :
\bigskip
\noindent
{\bf Proposition 1.1.} {\it The minimal polynomials and characteristic polynomials of  isometries of  quadratic spaces
are $\epsilon$--symmetric, where $\epsilon$ is the constant term of the polynomial.}
\medskip
\noindent
{\bf Proof.} Let $(V,q)$ be a quadratic space, and let $t : V \to V$ be
an isometry of $q$. By definition, we have  $q(tx,y) = q(x,t^{-1}y)$ for all $x,y \in V$. 
This implies that
for any polynomial $p \in k[X]$, we have $ q(p(t)x,y) = q(x,p(t^{-1})y)$  for all $x,y \in V$.
Let $f \in k[X]$ be the minimal polynomial of $t$. 
Applying the above equality to $p=f$, we see that the endomorphism $t^{{\rm deg}(f)} f(t^{-1})$ annihilates $V$.
As $f$ is the minimal polynomial of $t$, this implies that $f$ divides $X^{{\rm deg}(f)} f(X^{-1})$,
therefore we have $f(X) = \epsilon X^{{\rm deg}(f)} f(X^{-1})$ for some $\epsilon = \pm 1$.
On the other hand, the coefficient of the leading term of $X^{{\rm deg}(f)} f(X^{-1})$ is
equal to $f(0)$. Therefore $\epsilon = f(0)$. The statement
concerning the characteristic polynomial follows from a straightforward computation,
see for instance [Le 69], Lemma 7, (a).
\bigskip
If $f \in k[X]$ is a monic polynomial such that $f(0) \not = 0$, set
 $$f^*(X) = { 1 \over {f(0)}} X^{{\rm deg}(f)} f(X^{-1}).$$
 Note that $f^*$ is also monic, that $f^*(0) \not = 0$, and that $f^{**} = f$.

\bigskip
\noindent
{\bf Definition 1.2.} Let $f \in k[X]$ be a monic, $\epsilon$--symmetric polynomial with $\epsilon = \pm 1$. We say that $f$ is of 
\medskip
\noindent
$\bullet$ {\it type} 0 if $f$ is a product of powers of $X-1$ and of $X+1$;
\medskip
\noindent
$\bullet$ {\it type} 1  if $f$ is a product of powers of monic, symmetric, irreducible polynomials in
$k[X]$ of even degree;
\medskip
\noindent
$\bullet$ {\it type} 2 if $f$ is a product of polynomials of the form $g g^*$, where $g \in k[X]$ is
monic, irreducible, and $g \not = \pm  g^*$. 

\bigskip
\noindent
{\bf Proposition 1.3.}  \ {\it Every monic, $\epsilon$--symmetric polynomial  $F \in k[X]$ is a product of polynomials
of type 0, 1 and 2.}
\medskip
\noindent
{\bf Proof.} Let $f \in k[X]$ be a monic, irreducible factor of $F$. If $f \not = \pm  f^*$, then $f^*$ also
divides $F$, hence we get a factor of type 2. Suppose that $f = \pm f^*$. It suffices to show that
if $f(X) \not = X - 1, X + 1$, then ${\rm deg}(f)$ is even and $\epsilon = 1$. We have 
$f(X) = \epsilon X^{{\rm deg} f} f(X^{-1})$ for some $\epsilon = \pm 1$. If $\epsilon = -1$, then 
$f(1) = 0$, hence $f$ is divisible by $X-1$, and this is impossible as $f$ is supposed to be
irreducible and $f(X) \not = X - 1$. Hence $\epsilon = 1$. If  ${\rm deg}(f)$ is odd, then this implies that $f(-1) = 0$, which contradicts the assumption that $f$ is irreducible and $f(X) \not = X+1$.
Therefore ${\rm deg}(f)$
is even.

\bigskip
We say that a monic, $\epsilon$--symmetric polynomial is {\it hyperbolic} if all its components of
type 0 and 1 are of the form $f^e$ with $e$ even.

\bigskip

{\bf \S 2. Self--dual torsion modules}

\bigskip
Let $V$ be a finite dimensional $k$--vector space, and let $t  : V \to V$ be an endomorphism.
Then $V$ has a structure of torsion $k[X]$--module obtained by setting $X.v = t(v)$ for
all $v \in V$. 
Let us denote by $M(t)$ the
torsion $k[X]$--module associated to the endomorphism $t$. The module $M(t)$ will
be called the {\it module of the endomorphism} $t$. 
\medskip
Any torsion $k[X]$--module is isomorphic to a direct sum of modules of the
form $[k[X]/(f)]^m$ for some $f \in k[X]$ and $m \in {\bf N}$. 
If $M$ is a torsion $k[X]$--module,
set $F_M = \prod f^m$ for all $f \in k[X]$ and
$m \in {\bf N}$ as above. We call $F_M$ the {\it characteristic polynomial} of M. Note that
if $M = M(t)$ for some endomorphism $t$, then $F_M$ is the
characteristic polynomial of $t$. 
\medskip
A torsion $k[X]$--module is said to be of {\it type i}, for $i = 0,1,2$, if 
$M$ is a direct sum of modules of the form [$k[X]/(f)]^m$ where $f \in k[X]$ is of type $i$
and $m \in {\bf N}$.
It is said to be {\it self--dual} if it is a direct sum of modules of type 0, 1 and 2.
\medskip
From now on, {\it module} will mean a self--dual torsion $k[X]$--module that is
finite dimensional as a $k$--vector space.
\medskip
A module is said to be {\it hyperbolic} if
all its components of type 0 and type 1 are of the form $[k[X]/(f^e)]^n$ with $e$ even. 
We will see that any quadratic space having an isometry with hyperbolic module is 
hyperbolic.

\bigskip
{\bf  \S 3. Primary decomposition and transfer}

\bigskip
The aim of this section is to recall some results of Milnor [M 69].
\bigskip
Let $(V,q)$ be a quadratic space of dimension $2n$, let $t$ be an isometry of $q$ and
let $F$ be the characteristic polynomial of $t$. For each monic, irreducible factor $f$ of $F$,
set $$V_f = \{ v \in V \ | \ f^i(t)(v) = 0 \ {\rm for \ some} \  i \in {\bf N} \}.$$

Let  $U$ and $W$ be
two subspaces of $V$. We say that $U$ and $W$ are {\it orthogonal} to each other
if $q(u,w) = 0$ for all $u \in U$ and $w \in W$. 
 We say
that $(V,q)$ is {\it hyperbolic} if $V$ has a self--orthogonal subspace of dimension $n$.
\bigskip
\noindent
{\bf Proposition 3.1.} {\it 
Let $f$ and $g$ be two monic, irreducible factors of $F$.
If $f \not = g^*$, then $V_f$ and $V_g$ are orthogonal to each other.}
\medskip
\noindent
{\bf Proof.} See Milnor [M 69], Lemma 3.1. 
\bigskip
\noindent
{\bf Corollary 3.2.} {\it If $f$ is not symmetric, then $(V_f \oplus V_{f^*},q)$ is hyperbolic.}
\medskip
\noindent
{\bf Proof.} See [M 69], \S 3, Case 3.
\bigskip
\noindent
{\bf Proposition 3.3.} {\it We have the following orthogonal decomposition
$$(V,q) \simeq \bigoplus (V_f,q) \oplus H$$
where the sum is taken over all distinct monic, symmetric and irreducible factors of F, 
and where $H$
is a hyperbolic space. The orthogonal factors are stable by the isometry.}
\medskip
\noindent
{\bf Proof.} This follows from prop. 3.1 and cor. 3.2. 
\bigskip
\noindent
{\bf Proposition 3.4.} {\it Let $f \in k[X]$ be monic, symmetric and irreducible. The quadratic
space $(V_f,q)$ decomposes as an orthogonal sum of factors with  modules
of the form  $$ [k[X]/(f^e)]^m$$
for some integers $e$ and $m$. If $e$ is even, then this orthogonal factor is hyperbolic.
If $e$ is odd, then it is the orthogonal sum of a hyperbolic space and of a quadratic
space having an isometry with module $[k[X]/(f)]^m$. Conversely, any quadratic space having an isometry with module $ [k[X]/(f)]^m$ and
any $e \in {\bf N}$ arises from a quadratic space with isometry of module
$ [k[X]/(f^e)]^m$ in this way.}
\medskip
\noindent
{\bf Proof.} See [M 69],  th. 3.2,  3.3 and 3.4.
\bigskip
Recall that {\it module} means a self--dual torsion $k[X]$--module which is finite
dimensional as a $k$--vector space, and see def. 2.2 for the definition of a hyperbolic
module. 
\bigskip
\noindent
{\bf Corollary  3.5.} {\it A quadratic space having an isometry with hyperbolic
module is hyperbolic.}
\medskip
\noindent
{\bf Proof.} This follows from prop. 3.3 and 3.4. 
\bigskip
\noindent
{\bf Proposition 3.6.} {\it Let $f\in k[X]$ be a monic, symmetric and irreducible polynomial, and set
$K = k[X]/(f)$. Then sending $X$ to $X^{-1}$ induces a $k$--linear involution of $K$ denoted
by $x \mapsto \overline x$. Let $\ell : K \to k$ be a non--trivial $k$--linear map such that 
$\ell (x) = \ell (\overline x)$ for all $x \in K$. 
Then 
for every quadratic space $(V,q)$ over $k$ and every isometry 
having minimal polynomial $f$, there
exists a non--degenerate hermitian form $(V,h)$ over $K$ such that for all $x,y \in V$, we have
$$q(x,y) = \ell (h(x,y)).$$ 
Conversely, if $V$ is a finite dimensional vector space over $K$ and if $h : V \to V$ is
a non--degenerate hermitian form, then setting 
$$q(x,y) = \ell (h(x,y))$$
for all $x,y \in V$ we obtain a  quadratic space $(V,q)$ over $k$ together with an isometry
with minimal polynomial $f$. }
\bigskip
\noindent
{\bf Proof.} This is proved in [M 69], Lemma 1.1 and Lemma 1.2, in the case where $f$
is separable and 
$\ell  = {\rm Tr}_{K/k}$ is the trace of the extension $K/k$. The proof is the same for any non--trivial linear map $\ell$ having
the property that $\ell (x) = \ell  (\overline x)$ for all $x \in K$, as pointed out in [M 69], Remark 1.4.

\bigskip
\noindent
{\bf Corollary 3.7.} {\it Let $M$ be a module. Then there exists a quadratic space $q$ having
an isometry $t$ such that $M(t) \simeq M$.}
\medskip
\noindent
{\bf Proof.} This follows from prop. 3.3, 3.4 and 3.6. 

\bigskip
This implies the following well--known fact :
\bigskip
\noindent
{\bf Corollary 3.8.} {\it Let $F \in k[X]$ be a monic, $\epsilon$--symmetric polynomial. Then there
exists a quadratic space having an isometry with characteristic polynomial $F$.}
\medskip
Note that a new proof of this result, based on Bezoutians, is given in [JRV 12], th. 4.1. Integral analogs of
this question are investigated in [JRV 12], \S 3, as well as in [B 84], [BMar 94] and [B 99].

\bigskip
{\bf \S 4. Isometries with a given module}
\bigskip
We keep the notation of the previous sections. In particular,  {\it module}  means  a
self--dual torsion $k[X]$--module which is a finite dimensional $k$--vector space. 
\bigskip
The aim of this paper  is to investigate the following question
\bigskip
\noindent
{\bf Question.}  Which quadratic spaces admit an isometry of given module ?

\bigskip
This is a generalization of Milnor's question quoted in the introduction. 
Let us fix some notation. For
any module $M$, we have  $$M = M^0 \oplus M^1 \oplus M^2,$$
with $M^i$ of type $i$. Note that ${\rm dim}(M^2)$ is even, and let 
$2m_2 = {\rm dim}(M^2)$. Let us write
$M^0 = M_+ \oplus M_-$, with $M_+ = [k[X]/(X+1)^{e_+}]^{m_+}$ and 
$M_- = [k[X]/(X-1)^{e_-}]^{m_-}$, for some $e_+, e_-, m_+, m_- \in {\bf N}$.  Set
$n_+ = e_+ m_+$ and $n_- = e_- m_-$. Note that ${\rm dim}(M^0) = n_+ + n_-$. 

\bigskip

Let us first prove  that it is sufficient to consider
{\it semi--simple} modules. For any module $M$, let us denote by  ${\rm rad( M)}$ its radical,
and set $\overline M = M/{\rm rad}(M)$. 
Any module $M$ is the direct sum of modules of the form $M_{f,e} = [k[X]/(f^e)]^n$ for
some $f \in k[X]$ with $f$ irreducible and $e,n \in {\bf N}$. Let $M_{\rm odd}$ be
the direct sum of the modules $M_{f,e}$ with $e$ odd and $f$ symmetric. 
\bigskip
\noindent
{\bf Proposition 4.1.} {\it Let $q$ be a quadratic space and let $M$ be a module. Then
$q$ has an isometry with module $M$ if and only if $q$ is isomorphic to the orthogonal
sum of a quadratic space $\overline q$ with module $\overline M_{\rm odd}$ and of a
hyperbolic space.}

\medskip
\noindent
{\bf Proof.} Suppose that $q$ has an isometry with module $M$. Then by 3.3, 3.4 and 3.5, the
 quadratic space $q$ is isomorphic to the orthogonal sum of a quadratic space $q_{\rm odd}$ having
 an isometry with module $M_{\rm odd}$ and of a hyperbolic space. Moreover, by 3.4
 the quadratic space $q_{\rm odd}$ is isomorphic to the orthogonal sum of quadratic
 spaces $q_{f,e}$ having isometries with modules $M_{f,e}$.  Further, 3.4 also
 implies that $q_{f,e} \simeq \overline q_{f,e} \oplus H_{f,e}$, where $H_{f,e}$ is a hyperbolic
 space and $ \overline q_{f,e}$ has an isometry with  module $\overline M_{f,e}$.
 Note that $\overline M_{\rm odd}$ is the direct sum of the modules $\overline M_{f,e}$
 for $e$ odd and $f$ symmetric. Let $\overline q$ be the orthogonal sum of the
 quadratic spaces $ \overline q_{f,e}$. Then $\overline q$ has an isometry with
 module $\overline M_{\rm odd}$, and $q$ is the orthogonal sum of $\overline
 q$ and of a hyperbolic space.

 \medskip
 Conversely, suppose that $q \simeq \overline q \oplus H$, where $\overline q$ is a
 quadratic space having an isometry with module $\overline M$, and where $H$ is
 a hyperbolic space.  Then $\overline q$ is the orthogonal sum of quadratic
 spaces $ \overline q_{f,e}$ having isometries with modules $\overline M_{f,e}$.
 By 3.4 we get quadratic spaces $q_{f,e} \simeq \overline q_{f,e} \oplus H_{f,e}$
 having isometries with modules $M_{f,e}$, and $q$ is the orthogonal sum
 of the spaces $q_{f,e}$ and of a hyperbolic space. Hence $q$ has an
 isometry with module $M$.

 \bigskip
 
Recall that
the {\it Witt index} of a quadratic space $q$ is the number of hyperbolic planes contained in
the Witt decomposition of $q$. 

\bigskip
\noindent
{\bf Lemma 4.2.} {\it If $q$ is a quadratic space having an isometry with module $M$, then the Witt index
of $q$ is $\ge m_2$. }

\bigskip
\noindent
{\bf Proof.} Indeed, by
3.5, any quadratic space having an isometry with a module of type 2 is
hyperbolic. 

\bigskip
For the remainder of this section, let us assume that $M$ is a {\it semi--simple} module. Then the converse also holds if $M^1 = 0$ :

\bigskip
\noindent
{\bf Proposition 4.3.} {\it Let $q$ be a quadratic space such that ${\rm dim}(q) = {\rm dim}(M)$. Suppose that 
$M^1 = 0$. Then  $q$ has an isometry with module $M$ if and only if the
Witt index of $q$ is at least equal to $m_2$.}
\bigskip
\noindent
{\bf Proof.} We already know that if $q$ has an isometry with module $M$, then
Witt index of $q$ is at least  $m_2$. 
Conversely, suppose that the Witt index of $q$ is at least  $ m_2$,
and write $q$ as the orthogonal sum
of a quadratic space $(V_0,q_0)$ with a hyperbolic form of dimension 
$2m_2$. Then  ${\rm dim}(V_0) = {\rm dim}(M^0) =  n_+ + n_-$. Let us decompose
$(V_0,q_0)$ as an orthogonal sum of $(V_+,q_+)$ and $(V_-,q_-)$ with ${\rm dim}(V_+) = n_+$,
 ${\rm dim}(V_-) = n_-$, and let $t : V_0 \to V_0$ be defined by $t(x) = -x$ if $x \in V_+$
and $t(x) =x$ if $x \in V_-$. Then $t$ is an isometry of $q_0$, hence  we obtain
an isometry of $q$ with module $M$. 
\bigskip
This proposition has some useful consequences. In order to state the first one, we need
the notion of $u$--invariant of a field. 
Recall that a quadratic space $(V,q)$ is said to be {\it isotropic} if there
exists a non--zero $v \in V$ with $q(v,v) = 0$, and {\it anisotropic} otherwise. The
{\it $u$--invariant} of $k$, denoted by $u(k)$, is the largest dimension of an anisotropic
quadratic form over $k$. 
\bigskip
\noindent
{\bf Corollary 4.4.} {\it Let $q$ be a quadratic space with ${\rm dim}(q) = {\rm dim}(M)$.  Suppose that $u(k) \le {\rm dim}(M^0)$. Then 
$q$ has an isometry with module $M$.}
\medskip
\noindent
{\bf Proof.} Let $q_1$ be a quadratic space having an isometry with module $M^1$ (cf. prop. 3.7).
As $u(k) \le  {\rm dim}(M^0)$, we have  $q \oplus (-q_1) \simeq q_0 \oplus H$, where $H$ is
hyperbolic and 
$q_0$ is a quadratic space with ${\rm dim}(q_0)
=  {\rm dim}(M^0)$. 
Let $q_2$ be the hyperbolic space of dimension $2m_2$. 
Then the quadratic space $q_0 \oplus q_2$ has dimension ${\rm dim}(M^0 \oplus
M^2)$ and Witt index $ \ge m_2$. Therefore by prop. 4.3
the quadratic space $q_0 \oplus q_2$ has an isometry with module $M^0 \oplus
M^2$; hence  $q_0 \oplus q_1 \oplus q_2$ has an isometry with
module $M$. We have $q \oplus q_2 \oplus q_1 \oplus (-q_1) \simeq q_0 \oplus q_1 \oplus
q_2 \oplus H$. Since  $q_1 \oplus (-q_1)$ and $q_2$ are hyperbolic, and ${\rm dim}(q) = {\rm dim}(M) =
{\rm dim}(q_0 \oplus q_1 \oplus q_2)$, Witt cancellation implies that
$q \simeq q_0 \oplus q_1 \oplus q_2$. Therefore $q$ has an isometry with
module $M$. 
\bigskip
Let us recall that two quadratic spaces $q$ and $q'$ are {\it Witt--equivalent} if 
there exist hyperbolic spaces $H$ and $H'$ such that $q \oplus H \simeq q' \oplus H'$.
\bigskip
\noindent
{\bf Corollary 4.5.} {\it Suppose that $M^1 = 0$, and let $q$ be a quadratic space with ${\rm dim}(q) 
= {\rm dim}(M^0)$. Then any quadratic space of dimension equal to ${\rm dim}(M)$ and 
Witt--equivalent to $q$ has an isometry with module $M$.}
\medskip
\noindent
{\bf Proof.} Indeed, as $M^1 = 0$ we have $M = M^0 \oplus M^2$, hence 
${\rm dim}(M) = {\rm dim}(M^0) + 2m_2.$ Let $q'$ be a quadratic
space with ${\rm dim}(q) = {\rm dim}(M)$ and Witt--equivalent to $q$. Then 
the Witt index of $q'$ is at least $m_2$, hence by prop. 4.3
the quadratic space 
$q'$ has an isometry with module $M$. 
\bigskip
The next corollary will be used several times in the sequel.  
\bigskip
\noindent
{\bf Corollary 4.6.} {\it Suppose that $M^0 \not = 0$, and let
$d \in k^*$.  Then there exists a quadratic space $q$ having an isometry with module $M$
and determinant $d$.}

\medskip
\noindent
{\bf Proof.} We have $M = M^0 \oplus M^1 \oplus M^2$. If $M^1 \not = 0$, let $q_1$ be a 
quadratic space having an isometry with module $M^1$ (cf. prop. 3.7), and let $d_1 = {\rm det}(q_1)$.
If $M^1 = 0$, set $d_1 = 1$. 

\medskip
Since $M^0 \not = 0$, there exists a quadratic
space $q_0$ with determinant $d d_1(-1)^{m_2}$  and ${\rm dim}(q_0) = {\rm dim}(M^0)$.
Let $H$ be the hyperbolic form of dimension $2m_2 = {\rm dim}(M^2)$, and set 
$q_2 = q_0 \oplus H$. Then ${\rm dim}(q_2) = {\rm dim}(M^0 \oplus M^2)$, and
${\rm det}(q_2) = d d_1$. Moreover, $q_2$ is Witt equivalent
to $q_0$. Therefore cor. 4.5 implies that the quadratic space $q_0$ has an isometry with module $M^0 \oplus M^2$.
Set $q = q_1 \oplus q_2$. Then ${\rm det}(q) = d \in k^*/k^{*2}$, and
$q$ has an isometry with module $M$. 
\bigskip 
\bigskip
{\bf \S 5. Determinants and values of the characteristic polynomial}

\bigskip We have a relationship between the determinant of the quadratic
space and the values of the characteristic polynomials of its isometries, as follows  : 
\bigskip
\noindent
{\bf Proposition 5.1.} {\it Let $(V,q)$ be a quadratic space, 
and let $F \in k[X]$ be
the characteristic polynomial of an isometry $t$ of $q$. Then we have

$${\rm det}(q)F(1) F(-1) \in k^2.$$}

\medskip
\noindent
{\bf Proof.} Let us define $q' : V \times V \to k$ by $q'(x,y) = q(x,(t-t^{-1})(y))$. Then
$q'$ is skew--symmetric, hence ${\rm det}(q') \in k^2$. On the other hand, we have
$${\rm det}(q') =  {\rm det}(q) {\rm det}(t) {\rm det}(t+1) {\rm det}(t-1) = {\rm det}(q) {\rm det}(t)F(1) F(-1).$$ If $F(1) F(-1) = 0$, then the
statement is clear, so we can assume that $F(1) F(-1) \not = 0$. It is easy to see that
$F(1) \not = 0$ implies that $F(X) = X^{{\rm deg}( F) } F(X^{-1})$  and that $F(-1) \not = 0$
implies that ${\rm deg}(F)$ is even. Hence ${\rm det}(t) = 1$, and we have
${\rm det}(q)F(1) F(-1) \in k^2$, as stated. 

\bigskip

The
following corollary is well--known  (see  for instance [Le 69], Lemma 7, (c), or  the appendix of [GM 02]) :
\bigskip
\noindent
{\bf Corollary 5.2.} {\it Let $q$ be an quadratic space, 
and let $F \in k[X]$ be
the characteristic polynomial of an isometry  of $q$. Suppose that $F(1)F(-1) \not = 0$. Then

$${\rm det}(q) = F(1) F(-1) \in k^*/k^{*2}.$$}
\medskip
\noindent
{\bf Proof.} This is an immediate consequence of 5.1.

\bigskip
The following lemma will be useful in the next sections.  \bigskip
\noindent
{\bf Lemma 5.3.} {\it Let  $M$ be a semi--simple module, and let $d \in k^*$ such that 
$d F_M(1) F_M(-1) \in k^2.$ Then there exists a
quadratic space $q$ of determinant $d$ having an isometry with module $M$. }
\medskip
\noindent
{\bf Proof.} Suppose first that 
$F_M(1)F_M(-1) \not = 0$; then the hypothesis implies that 
$d = F_M(1) F_M(-1) \in k^*/k^{*2}$. Let $q$ be any
quadratic space with module $M$ (cf. prop. 3.7). By cor. 5.2 we have 
${\rm det}(q) = F_M(1) F_M(-1)$ in $k^*/k^{*2}$, hence 
${\rm det}(q) = d$ in $k^*/k^{*2}$.  Suppose now that $F_M(1)F_M(-1) = 0$, and note that
this implies that $M^0 \not = 0$.  By cor. 4.6, there exists 
a quadratic space $q'$ with determinant $d$ having an isometry of module $M$,
and this completes the proof of the lemma.

\bigskip
{\bf \S 6. Fields with $I(k)^2 = 0$.}
\bigskip
We keep the notation introduced in \S 4. In particular, {\it module} means a 
self--dual torsion $k[X]$--module that is a finite dimensional $k$--vector space.
Recall that by prop. 4.1 it is sufficient to consider the case of {\it semi--simple}
modules. 
Let $W(k)$ be the Witt ring of $k$, and let $ I(k)$ be the fundamental ideal
of $W(k)$. 
Let $q$ be a quadratic space, let $M$ be
a semi--simple module such that 
${\rm dim}(q) = {\rm dim}(M)$.

\bigskip
\noindent
{\bf Proposition 6.1.} {\it Suppose that $I(k)^2 = 0$. Then the quadratic space $q$ has an isometry with module $M$ if
and only if 
$${\rm det}(q)F_M(1) F_M(-1) \in k^2.$$}

\medskip
\noindent
{\bf Proof.} The condition is necessary by prop. 5.1. Let us show that it is sufficient. 
 By lemma 5.3 there
exists a quadratic space $q'$ having an isometry with module $M$ and such that
${\rm det}(q') = {\rm det}(q)$. 
Then $q$ and $q'$ have the same dimension and determinant. As $I(k)^2 = 0$, this implies that they
are isomorphic, therefore $q$ has an isometry with module $M$.

\bigskip
\noindent
{\bf Corollary 6.2.}   {\it Suppose that $I(k)^2 = 0$, and that $F_M(1) F_M(-1) \not = 0$.  Then the quadratic space $q$ has an isometry with module $M$ if
and only if 
$${\rm det}(q) = F_M(1) F_M(-1) \in k^*/k^2.$$}

\medskip
\noindent
{\bf Proof.} This follows from 6.1.

\bigskip
Let $k_s$ be a separable closure of $k$, and set $\Gamma_k = {\rm Gal}(k_s/k)$. We
say that the $2$--cohomological dimension of $k$, denoted by ${\rm cd}_2(k)$, is at most  1 if $H^r
(\Gamma_k,A) = 0$ for all finite $2$--primary $\Gamma_k$--modules $A$ and
for all $r > 1$. 

\bigskip
\noindent
{\bf Corollary 6.3.} {\it Suppose that ${\rm cd}_2(k) \le 1$. Then the quadratic space $q$ has an isometry with module $M$ if
and only if $${\rm det}(q)F_M(1) F_M(-1) \in k^2.$$
If moreover $F_M(1) F_M(-1) \not = 0$, then $q$ has an isometry with module $M$ if
and only if $${\rm det}(q) = F_M(1) F_M(-1) \in k^*/k^2.$$}

\medskip
\noindent
{\bf Proof.}  It is well--known that if ${\rm cd}_2(k) \le 1$, then $I(k)^2 = 0$, hence
the corollary follows from th. 6.1. and cor. 6.2.

\bigskip
{\bf \S 7. Local fields}
\bigskip
We keep the notation of \S 4. 
In particular, {\it module} means a 
self--dual torsion $k[X]$--module that is a finite dimensional $k$--vector space.
For any module $M$, we have $M = M^0 \oplus M^1 \oplus M^2$, where $M^i$ is of type $i$. 
Let us suppose that $M$ is {\it semi--simple} (this is possible by prop. 4.1). Note
that if $M^1 = 0$, then a quadratic space has an isometry with module $M$ if
and only if its Witt index is $\ge m_2$, where $2m_2 = {\rm dim}(M^2)$ (cf. prop. 4.3).
Therefore from now on we can restrict ourselves to modules  $M$  with $M^1 \not = 0$. 

\bigskip
Suppose that $k$ is a local field. Let $q$ be a quadratic space, and let $M$ be a module with $M_1 \not = 0$. Suppose that
${\rm dim}(q) = {\rm dim}(M)$. 

\bigskip
\noindent
{\bf Theorem 7.1.} {\it The quadratic space $q$ has an isometry with module $M$ if and 
only if $${\rm det}(q)F_M(1) F_M(-1) \in k^2.$$}

The proof of th. 7.1 uses the following result of Milnor. 
 Let $K$ be an extension of $k$ of finite degree endowed
with a non--trivial $k$--linear involution $x \mapsto \overline x$. Let $\ell : K \to k$ be
a non--trivial linear form such that $\ell (x) = \ell (\overline x)$ for all $x \in K$.
For any non--degenerate hermitian form $h : V \times V \to K$, let us denote by
$q_h : V \times V \to k$ the quadratic space defined by 
$q_h(x,y) = \ell (h(x,y))$ for all $x,y \in V$.  We have
\bigskip
\noindent
{\bf Theorem 7.2.}  {\it If the hermitian spaces $h$ and $h'$ have the
same dimension but different determinants, then the quadratic spaces $q_h$ and
$q_{h'}$ have the same dimension and determinant but different Hasse invariants.}
\medskip
\noindent
{\bf Proof.} See Milnor, [M 69], th. 2.7.
\bigskip
\noindent
{\bf Proof of th. 7.1} If $q$ has an isometry with  module $M$, then by prop. 5.1
we have ${\rm det}(q)F_M(1) F_M(-1) \in k^2.$ Conversely, let us
suppose that $${\rm det}(q)F_M(1) F_M(-1) \in k^2.$$
By lemma 5.3, there exists a quadratic space $q'$ having an
isometry with module $M$ such that ${\rm det}(q') = {\rm det}(q)$. It is well--known
that two quadratic spaces over a local field are isomorphic if and only
if they have the same dimension, determinant and Hasse--Witt invariant. Therefore
if the Hasse--Witt
invariants of $q$ and $q'$ are equal, then $q \simeq q'$, hence we are finished.
\medskip
Suppose that this is not the case. As $M^1 \not = 0$, there exists a monic, symmetric,
irreducible polynomial $f \in k[X]$ of even degree such that for some $n \in {\bf N}$
and for some odd integer $e$, the module $M_f = [k[X]/(f^e)]^n$ is a direct summand
of $M^1$. Set $M = M_f \oplus \tilde M$. By prop. 3.3 and 3.4, we have an orthogonal decomposition 
$q' \simeq q_f \oplus \tilde q$, where $q_f$ has an isometry with module $M_f$ and $\tilde q$ has
an isometry with module $\tilde M$. 
\medskip
Set $K = k[X]/(f)$, and let us consider the $k$--linear involution of $K$ induced by
$X \mapsto X^{-1}$. Let $E$ be the fixed field of this involution.  Set $V = K^n$.
Then by prop. 3.4 and
3.6, there exists a hermitian form $h : V \times V \to K$ such that the quadratic space $q_f$ has an orthogonal decomposition $q_f \simeq q_h \oplus H$,
where $H$ is a hyperbolic space and where $q_h : V \times V \to k$ is defined by
$$q_h(x,y) = \ell  (h(x,y)).$$
Let $\alpha_1,\dots,\alpha_n \in E^*$ be such that $h \simeq <\alpha_1,\dots,\alpha_n>$. Let
us denote by $N_{K/E} : K \to E$ the norm map, 
and let $\alpha \in E^*$ be such that $\alpha \not \in N_{K/E}(K^*)$. Let $h' : V \times V \to K$
be the hermitian form defined by $h ' =  <\alpha \alpha_1,\dots,\alpha_n>$. 
Let us define $q_{h'} : V \times V \to k$  by
$$q_{h'}(x,y) = \ell  (h'(x,y)).$$
Then
$h$ and $h'$ have the same dimension but different determinants. Therefore by th. 7.2, 
the quadratic forms $q_h$ and $q_{h'}$ have the same dimension and determinant
but different Hasse--Witt invariants. 
\medskip
Set $q'_f = q_{h'} \oplus H$. By prop. 3.4, the quadratic space $q'_f$ has an isometry
with module $M_f$. Let $q'' = q'_f \oplus \tilde q$. Then $q''$ has an isometry with
module $M$, and the quadratic spaces $q$ and $q''$ have equal dimension,
determinant and Hasse--Witt invariants. Therefore $q \simeq q''$, hence $q$ has
an isometry with module $M$ as claimed. 
\bigskip
Recall that we are assuming that $M^1 \not = 0$. The following corollary shows
that if in addition $M^0 \not = 0$, then any quadratic form of dimension ${\rm dim}(M)$
has an isometry with module $M$. 
\bigskip
\noindent
{\bf Corollary 7.3.} {\it Suppose that $M^0 \not = 0$. Then any quadratic space of
dimension ${\rm dim}(M)$ has an isometry with module $M$. }
\medskip
\noindent
{\bf Proof.} Let $q$ be a quadratic space with ${\rm dim}(q) =
{\rm dim}(M)$. As $M^0 \not = 0$, we have $F_M(1)F_M(-1) = 0$, therefore the
condition ${\rm det}(q)F_M(1) F_M(-1) \in k^2$ holds independently
of the value of ${\rm det}(q)$. Hence by th. 7.1 the quadratic form $q$ has
an isometry with module $M$. 

\bigskip
\noindent
{\bf Corollary 7.4.}   {\it Suppose that  $F_M(1) F_M(-1) \not = 0$.  Then the quadratic space $q$ has an isometry with module $M$ if
and only if 
$${\rm det}(q) = F_M(1) F_M(-1) \in k^*/k^2.$$}

\medskip
\noindent
{\bf Proof.} This is a consequence of 7.1.

\bigskip
{\bf \S 8. The field of real numbers}
\bigskip
In this section the ground field $k$ is the field of real numbers ${\bf R}$. Let $q$ be
a quadratic space over ${\bf R}$. It is well--known that $q$ is isomorphic to 
$$X_1^2 +  \dots + X_r^2 - X_{r+1}^2 - \dots - X_{r + s}^2$$
for some natural numbers $r$ and $s$. These are uniquely determined by $q$, and we
have $r + s = {\rm dim}(V)$. The couple $(r,s)$ is called the {\it signature} of $q$. 
\bigskip
Let $M$ be a module. Recall that 
$M = M^0 \oplus M^1 \oplus M^2$ with $M^i$ of type $i$. Let  
$F_M$ be the characteristic polynomial of $M$, and 
let $2 \sigma$ be the number of roots of $F_M$ off the unit circle. Note that 
${\rm dim}(M^2)  = 2 \sigma$.
\medskip
Let us introduce some notation. For any integers $n,m,n',m'$, we write $(n,m) \ge (n',m')$ if and
only if $n \ge n'$ and $m \ge m'$, and $(n,m) \equiv (n',m') 
\ \ ({\rm mod} \ 2)$ if and only if $n \equiv  n'
\ \ ({\rm mod} \ 2)$ and $m \equiv  m'
\ \ ({\rm mod} \ 2)$.
\medskip
We have seen in \S 4 that it suffices to consider semi--simple modules (cf. prop. 4.1). We first
give the criterion in the semi--simple case (see prop. 8.1 below), and then use prop. 4.1 to
treat the case of arbitrary modules.

\bigskip
\noindent {\bf Proposition 8.1.} {\it Assume that $M$ is semi--simple, and that  ${\rm dim}(q) = {\rm dim}(M)$. 
\medskip
\noindent
{\rm (a)}  Suppose that the quadratic space $q$ has an isometry
with module $M$. Then we have  $$(r,s) \ge (\sigma , \sigma).$$
If moreover $M^0 = 0$, then  $$( r,s) \equiv (\sigma , \sigma) 
\ \ ({\rm mod} \ 2).$$
\noindent
{\rm (b)}  
Conversely, suppose that  we have $$(r,s) \ge (\sigma, \sigma),$$
and that moreover if $M^0 = 0$, then
$$( r,s) \equiv (\sigma, \sigma) 
\ \ ({\rm mod} \ 2).$$
Then the quadratic space $q$ has an isometry with module $M$. } 
\medskip
\noindent
{\bf Proof.}  {\rm (a)}   Suppose that the quadratic space $q$ has an isometry
with module $M$. By prop. 3.3, we have an orthogonal decomposition
$$(V,q) \simeq \bigoplus (V_f,q_f) \oplus H$$
where the sum is taken over all distinct monic, symmetric and irreducible factors of $ F_M$,
and where $H$
is a hyperbolic space. Note that ${\rm dim}(H) = {\rm dim}(M^2) = 2 \sigma$. 
This implies that $(r,s) \ge (\sigma, \sigma).$
If $M^0 = 0$, then every irreducible and symmetric polynomial $f$
appearing in the above decomposition  is of degree two. Let $K_f = k[X]/(f)$. Then $V_f$
has a structure of $K_f$--vector space, and by prop. 3.6,
there exists a hermitian form $h_f: V_f \times V_f  \to K_f$ such that
$$q_i(x,y) = {\rm Tr}_{K/{\bf R}}(h_f(x,y))$$ for all $x,y \in V_{f,}$. Let $(u_f,v_f)$ be
the signature of $h_f$. Then the signature of $q_f$ is $(2u_f,2v_f)$, and this implies
that  $(r,s) \equiv (\sigma ,\sigma)
\ ({\rm mod} \ 2)$. 
\medskip
{\rm (b)}  Conversely, suppose that $r+ s = {\rm dim}(M) $ and that $(r,s) \ge (\sigma,\sigma)$.
Note that ${\rm dim}(M) - 2 \sigma  = {\rm dim}(M) - {\rm dim}(M_2)  \ge
{\rm dim}(M^0)$, therefore $r + s - 2 \sigma  \ge {\rm dim}(M^0) $. 
As ${\rm dim}(M^1)$ is even, we also have
${\rm dim}(M) = r + s \equiv {\rm dim}(M^0)  \ ({\rm mod} \ 2)$. 
Set $r' = r - \sigma$ and $s' = s - \sigma$. Then $r '+ s' \equiv 
r + s \equiv {\rm dim}(M^0)  \ ({\rm mod} \ 2)$, and $r'+ s' \ge {\rm dim}(M^0) $. 
Let us write $$r' = 2u + u_+  \  {\rm and}  \ s' =  2v + v_-$$ with $u,v, u_+, v_- \in {\bf N}$ such that
$u_+ + v_- = {\rm dim}(M^0)$.  Indeed, it is easy to see that this is possible if ${\rm dim}(M^0)  > 0$.
On the other hand, if  ${\rm dim}(M^0) = 0$ then  $M^0 = 0$,
hence by hypothesis $( r,s) \equiv (\sigma, \sigma) 
\ \ ({\rm mod} \ 2)$. This implies that
$r'$ and $s'$ are even. In this case, set $u = {r' \over 2}$ and $v = {s' \over 2}$. 
\medskip
Note that ${\rm dim}(M^1) = 2u + 2v$. 
Recall that $M^1$ is a direct sum of modules of the type $[k[X]/(f)]^{n_{f}}$ with
$f \in {\bf R}[X]$ symmetric, irreducible and ${\rm deg}(f) = 2$. 
Note that $u + v = \Sigma n_{f}$, where the sum is taken over all $f$
as above. 
Let
$u_{f},v_{f} \in {\bf N}$ such that $0 \le u_{f}, v_{f} \le n_{f}$, that $u_{f} + v_{f} = n_{f}$,
and that $$\Sigma u_{f} = u, \ \ \Sigma v_{f} = v,$$
the sums being taken over all the $f$ as above. 
\medskip
Set $K_f = k[X]/(f)$ and  $V_{f} = K_f^{n_{f}}$. Let $h_{f} : V_f \times V_f  \to K_f$
be a hermitian form of signature $(u_f,v_f)$, and let $q_{f} : V_f  \times V_f  \to {\bf R}$
be the quadratic space defined by $q_f (x,y) = {\rm Tr}_{K_f/{\bf R}}(h_{f}(x,y))$ for all
$x,y \in V_f $.  Then the signature of $q_f$ is $(2u_f, 2v_f)$. Let $q_1$ be the
orthogonal sum of the spaces $q_f$ for all $f$ as above. Then the signature of $q_1$ is
$(2u,2v)$. 
\medskip
Let $q_0$ be the quadratic space of signature $(u_+,v_-)$, and let $q_2$ be
the hyperbolic space of dimension $2 \sigma$. Let $q' = q_0 \oplus q_1 \oplus q_2$.
Then $q'$ has an isometry with module $M = M^0 \oplus M^1 \oplus M^2$.
\medskip
Note that the signature of $q'$ is $(u_+ + 2u + \sigma, v_- + 2v + \sigma) = (r,s)$. 
Therefore ${\rm sign}(q') = {\rm sign}(q)$, hence $q' \simeq q$. This implies that
$q$ has an isometry with module $M$.
\bigskip
\noindent
{\bf Corollary 8.2.} {\it Let $F \in {\bf R}[X]$ be a symmetric polynomial such that $F(1)F(-1) \not = 0$.
Then the quadratic space $q$ has a semi--simple  isometry with characteristic polynomial $F$ if and only if
$(r,s) \ge (\sigma, \sigma)$ and   $(r,s) \equiv  (\sigma, \sigma) \ ({\rm mod} \ 2)$.}
\medskip
\noindent
{\bf Proof.} Let $M$ be the semi--simple module with characteristic polynomial $F = F_M$. 
As
$F_M(1)F_M(-1) \not = 0$, we have  $M^0 = 0$, therefore the corollary follows
from prop. 8.1. 
\medskip
A special case of this corollary is proved by Gross and McMullen in  [GM 02], cor. 2.3.
\medskip
Prop. 8.1 and prop. 4.1 lead to  a criterion for arbitrary modules :
\medskip
\noindent
{\bf Corollary 8.3} {\it Suppose that  ${\rm dim}(q) = {\rm dim}(M)$,
and set $2 \tau = 
{\rm dim}(M) - 
{\rm dim}(\overline M^0_{\rm odd})$.
 \medskip
\noindent
{\rm (a)}  Suppose that the quadratic space $q$ has an isometry
with module $M$. Then we have  $$(r,s) \ge (\tau ,  \tau).$$
If moreover $\overline M^0_{\rm odd} = 0$, then  $$( r,s) \equiv   (\tau ,  \tau)
\ \ ({\rm mod} \ 2).$$
\noindent
{\rm (b)}  
Conversely, suppose that  we have $$(r,s) \ge  (\tau ,  \tau),$$
and that moreover if $\overline M^0_{\rm odd} = 0$, then
$$( r,s) \equiv  (\tau ,  \tau)
\ \ ({\rm mod} \ 2).$$
Then the quadratic space $q$ has an isometry with module $M$. } 
\medskip
\noindent
{\bf Proof.} {\rm (a)} As $q$ has an isometry with module $M$, by prop. 4.1 the quadratic space $q$ is
isomorphic to the orthogonal sum of a quadratic space $q'$ having an
isometry with module $\overline M_{\rm odd}$ and of a hyperbolic space $H$ of
dimension $2 \tau$.  The signature of $H$ is $(\tau,\tau)$,
hence  we have  $(r,s) \ge (\tau ,  \tau).$
Let $(r',s')$ be the signature of $q'$. 
Note that all the roots of the polynomial $F_{\overline M_{\rm odd}}$ are on the
unit circle. Therefore if  $\overline M^0_{\rm odd} = 0$, by prop. 8.1 we have
$(r',s') \equiv   (0,0)
\ \ ({\rm mod} \ 2)$, hence $( r,s) \equiv   (\tau,\tau)
\ \ ({\rm mod} \ 2).$

\medskip  {\rm (b)}  Since $( r,s) \equiv   (\tau,\tau)
\ \ ({\rm mod} \ 2)$, we have  $q \simeq q' \oplus H$, where $H$ is a hyperbolic space of
dimension $2 \tau$ and $q'$ is a quadratic space of dimension equal to ${\rm dim}(\overline M_{\rm odd})$. Let $(r',s')$ be the signature of $q'$.  If $\overline M^0_{\rm odd} = 0$, then by hypothesis we have $$( r,s) \equiv  (\tau ,  \tau)
\ \ ({\rm mod} \ 2).$$ Since the signature of $H$ is $(\tau,\tau)$, this implies that
$(r',s') \equiv   (0,0)
\ \ ({\rm mod} \ 2)$. Since all the roots of the polynomial $F_{\overline M_{\rm odd}}$ are on the
unit circle,  prop. 8.1 implies that the quadratic space $q'$ has an isometry
with module $\overline M_{\rm odd}$, and by prop. 4.1 this implies that $q$ has
an isometry with module $M$.

\bigskip
{\bf 9. Global fields - the case of an irreducible minimal polynomial}
\bigskip 
The aim of this section is to give an answer to Milnor's question
stated in the introduction in the case of global fields. 
Suppose that $k$ is a global field, let $q$ be a quadratic space, and let
$f \in k[X]$ be an irreducible and symmetric polynomial. We have the following Hasse principle :
\bigskip
\noindent
{\bf Theorem 9.1.} {\it The quadratic space $q$ has an isometry with minimal
polynomial $f$ if and only if such an isometry exists over every completion of $k$.}
\bigskip The case $f(X) = X+1$ is trivial, hence we  may assume that $f(1)f(-1) \not = 0$. Before proving th. 9.1, let us use 
the results of the previous two sections to obtain necessary and
sufficient conditions for an isometry to exist. Let $F$ be a power of $f$ 
such that ${\rm deg}(F) = {\rm dim}(V) = 2n$.
\bigskip
For every real place $v$ of $k$, let $(r_v,s_v)$ denote the signature of $q$ over $k_v$,
and let $\sigma_v$ be the number of roots of $F \in k_v[X]$ that are not on the unit circle. 
\bigskip
We say that the {\it signature condition} is satisfied for $q$ and $F$ if for every real place 
$v$ of $k$, we have $(r_v,s_v) \ge (\sigma_v,\sigma_v)$, and 
$(r_v,s_v) \equiv (\sigma_v,\sigma_v)
\ ({\rm mod} \ 2)$. 
\bigskip
We say that the {\it hyperbolicity condition} is satisfied for $q$ and $F$ if for all places
$v$ of $k$ such that $F \in k_v[X]$ is a hyperbolic polynomial, the quadratic form
$q_v$ over $k_v$ is hyperbolic. 
\bigskip
\noindent
{\bf Corollary 9.2.} {\it The quadratic space $q$ has an isometry with minimal
polynomial $f$ if and only if the signature condition and the hyperbolicity condition are satisfied, and 
$${\rm det}(q) = F(1)F(-1) \in k^*/k^{*2}.$$}
\medskip
\noindent
{\bf Proof.} The necessity of the conditions follows from prop. 8.2, 3.5  and 5.2.
Conversely,
suppose that the signature condition is satisfied and that ${\rm det}(q) = F(1)F(-1) \in k^*/k^{*2}.$
Then by prop. 8.2 and th. 7.4, the quadratic space $q$ has an isometry with minimal
polynomial $f$ over $k_v$ for every place $v$ of $k$. By th. 9.1, this implies that $q$
has an isometry with minimal polynomial $f$. 
\bigskip
The following reformulation of cor. 9.2 shows that it suffices to check a finite number
of conditions. Let $q$ and $F$ be as above, with ${\rm dim}(q) = {\rm deg}(F) = 2n$. 
Let $S$ be the set of places of $k$ at which the Hasse invariant of $q$ is not
equal to the Hasse invariant of the $2n$--dimensional hyperbolic space. Note
that $S$ is a finite set. 
\bigskip
\noindent
{\bf Corollary 9.3.} {\it The quadratic space $q$ has an isometry with minimal polynomial $f$
if and only if the following conditions are satisfied :
\medskip
{\rm (i)}  $F(1)F(-1) = {\rm det}(q) \in k^*/k^{*2}$;

{\rm (ii)} The signature condition holds;

{\rm (iii)} If $v \in S$, then $F \in k_v[X]$ is not hyperbolic.}

\medskip
\noindent
{\bf Proof.}  It suffices to prove that the conditions (i), {\rm (ii)} and (iii) imply the hyperbolicity condition.
Let $v$ be a place of $k$ such that $F \in k_v[X]$ is hyperbolic. Then there exists
a polynomial $G \in k_v[X]$ such that $F = G G^*$. Note that ${\rm deg}(G) = n$. 
We have $F(1) = G(1)^2$, and $F(-1) = (-1)^n G(-1)^2$. By condition (i), we have
$F(1)F(-1) = {\rm det}(q)$, hence $(-1)^n{\rm det}(q) = {\rm disc}(q) \in k_v^2$.
On the other hand, condition (iii) implies that $v \not \in S$, hence $q$ has the
same Hasse invariant  at $v$ as the $2n$--dimensional hyperbolic space.  
Therefore over $k_v$, the
quadratic space $q$ has the same dimension, discriminant and Hasse  invariant
as the $2n$--dimensional hyperbolic space. 
If $v$ is an infinite place, then by condition (ii) 
the signature of $q$  at $v$ coincides with the signature of the $2n$--dimensional hyperbolic space. Hence $q$ is hyperbolic over $k_v$, in
other words the hyperbolicity condition is satisfied.

\bigskip The following lemmas will be used in the proof of th. 9.1, and also in 
\S 10. 
\bigskip
Let $K = k[X]/(f)$, and let $^{\overline {\ }} : K \to K$ be the involution induced by
$X \mapsto X^{-1}$. Let $E$ be the fixed field of the involution. \bigskip
\noindent
{\bf Lemma 9.4.} {\it Let $v$ be a place of $k$. The following properties are equivalent :
\smallskip  
\noindent
{\rm (i)} Every place of $E$ above $v$
splits in $K$. 
\smallskip
\noindent
{\rm (ii)} The polynomial $f \in k_v[X]$ is hyperbolic.
\smallskip
\noindent
{\rm (iii)} For any $m \in {\bf N}$, the module $[k_v[X]/(f)]^m$ is hyperbolic.}
\bigskip
\noindent
{\bf Proof.} 
Let $w_1,\dots,w_r$ be the places of $E$ above $v$, and set
$E_i = E_{w_i}$, $K_i = K \otimes_E E_i$. Then 
$K_i$ is a field if $w_i$ is inert or ramified in $K$, a product
of two fields if $w_i$ is split in $K$, and $k_v[X]/(f) \simeq K_1 \times 
\dots \times K_r$.

\medskip 
\noindent
{\rm (i)} $\Longrightarrow$  {\rm (ii)} Since every $w_i$ splits
in $K$, all the $K_i$'s are products of two fields. This implies that 
$f = f_1f_1^* \dots f_r f_r^*$ with $f_i \in k_v[X]$ monic and irreducible and $f_i \not = f_i^*$
for all $i = 1,\dots,r$. Therefore $f \in k_v[X]$ is hyperbolic.
\smallskip 
\noindent
{\rm (ii)} $\Longleftrightarrow$ 
  {\rm (iii)} is clear.
 \smallskip
 \noindent
 {\rm (ii)} $\Longrightarrow$  {\rm (i)} Since $f \in k[X]$ is irreducible and $f \in k_v[X]$
 is hyperbolic, we have $f = f_1f_1^* \dots f_r f_r^*$ with $f_i \in k_v[X]$ monic and irreducible and $f_i \not = f_i^*$
for all $i = 1,\dots,r$. Therefore all the $K_i$'s are products of two fields, hence {\rm (i)} holds.

\bigskip
\noindent
{\bf Lemma 9.5.} {\it Let $v$ be a place of $k$ satisfying the equivalent conditions
of lemma 9.4. Then we have 
\smallskip
\noindent
{\rm (i)} Any quadratic space over $k_v$ having an isometry with minimal polynomial $f$  is hyperbolic.
\smallskip
\noindent
{\rm (ii)} For any $m \in {\bf N}$, every quadratic space over $k_v$ having an isometry with module $[k_v[X]/(f)]^m$ is
hyperbolic. }
\bigskip
\noindent
{\bf Proof.} Both assertions follow from lemma 9.4 and cor. 3.5.

\bigskip Recall that for any quadratic space $Q$, we denote by $w(Q)$ its Hasse
invariant. 
\bigskip
\noindent
{\bf Lemma 9.6.} {\it Let $m \in {\bf N}$, let $v$ be a finite place of $k$, and let 
$M =  [k_v[X]/(f)]^m$. Suppose that $M$
is not hyperbolic. Let $\epsilon \in \{ 0, 1 \}$. Then there exists a quadratic space $Q$
over $k_v$ such that $Q$ has an isometry with module $M$, and that
$w(Q) = \epsilon$.}

\bigskip
\noindent
{\bf Proof.} Note that $M$ is hyperbolic if and only if $f \in k_v[X]$ is
hyperbolic, in other words if it is a product of polynomials of type 2 over $k_v$. As
we are assuming that $M$ is not hyperbolic, the polynomial
$f \in k_v[X]$ has at least one irreducible factor of type 1. Hence we have
$f = f_1 f_2$ with $f_1, f_2 \in k_v[X]$ and $f_1$ irreducible, symmetric  of
even degree.
\medskip
Recall that  $K = k[X]/(f)$, that  $^{\overline {\ }} : K \to K$ is the involution induced by
$X \mapsto X^{-1}$ and that $E$ is the fixed field of this involution.
Set  $K_v = K \otimes _k k_v$. Then $K_v  \simeq K_1 \times K_2$ with $K_i = k_v[X]/(f_i)$,
the involution preserves $K_1$ and $K_2$, and we have $E = E_1 \times E_2$.
Note that $K_1$ is a field, and $E_1$ is the fixed field of the restriction of the
involution to $K_1$, hence $K_1/E_1$ is a quadratic extension. 
\medskip
We have $M   \simeq M_1 \oplus M_2$ with 
$M_1 \simeq K_1^m$ and $M_2 \simeq K_2^m$. Let $h : M \times M \to K_v$ be
the unit hermitian form. Then $h \simeq h_1 \oplus h_2$, where $h_i : M_i \times M_i \to K_v$,
with $i = 1,2$, is the restriction of $h$ to $M_i$.  Let $\ell : K_v  \to k_v$ be a non--zero
linear form such that $\ell (\overline x) = \ell (x)$ for all $x \in K_v$. For any hermitian form $H$, 
set $q_H(x,y)  = \ell (H(x,y))$. 
The quadratic space $q_h$ has an isometry with module $M$ by construction. If
$w(q_h) = \epsilon$, then we set $Q = q_h$ and the lemma is proved. 
\medskip
Suppose
that $w(q_h) \not = \epsilon$, and let $\alpha = {\rm det}(h_1)$; then $\alpha \in E_1^*$. 
Since $K_1/E_1$ is a quadratic extension, there exists $\beta \in E_1^*$ such that
$\beta \not \in {\rm N}_{K_1/E_1}(K_1^*)$. Let $h_1' : M_1 \times M_1 \to K_1$ be
a hermitian form of determinant $\alpha \beta$. Then $h_1$ and $h_1'$ have
same dimension and different determinants, hence by th. 7.2 the quadratic
spaces $q_{h_1}$ and $q_{h_1'}$ have equal dimension, determinant and different
Hasse invariants. Let $h' = h_1' \oplus h_2$, and set $Q = q_{h'}$. Then 
$Q \simeq q_{h_1'} \oplus q_{h_2}$, hence $w(Q) = \epsilon$. Since $Q$ has an
isometry with module $M$, this concludes the proof of the lemma.

 \bigskip
\noindent
{\bf Proof of th. 9.1} Let
 $F = f^m$ and  ${\rm deg}(f)  = 2d$, and recall that $n = md$. 
Let $K = k[X]/(f)$,  and let $^{\overline {\ }} : K \to K$ be the involution induced by
$X \mapsto X^{-1}$. Let $E$ be the fixed field of the involution. Let $\theta \in E^*$ be such that $K = E(\sqrt {\theta})$, and for any place $w$ of $E$,
 let $( \ , \ )_w$ denote the
Hilbert symbol at $E_w$. 
\bigskip
Let  $v$ be a real place of $k$. Then the signature $(r_v,s_v)$ of $q$ at $k_v$ 
satisfies  $(r_v,s_v) \ge (\sigma_v,\sigma_v)$ and 
$(r_v,s_v) \equiv (\sigma_v,\sigma_v)
\ ({\rm mod} \ 2)$.  In particular,  $r_v - \sigma_v $ is even. Set $r_v- \sigma_v = 2u_v$.
Then $s_v - \sigma_v = 2(n - \sigma_v - u_v)$, and $n - \sigma_v - u_v \ge 0$. 
Therefore we have $0 \le u_v \le n - \sigma_v$.
Let us denote by $2 \tau_v$ the number of roots of $f$ that are not on the unit circle. Then we have $\sigma_v = m \tau_v$.
Let us write $u_v = u_v^1 + \dots + u_v^m$ for some integers $u_v^i$ such that
$0 \le u_v^i \le d - \tau_v$.
\bigskip
Let $w_1,\dots,w_{d-\tau_v}$ be the real places of $E$ above $v$ that extend to
a complex place of $K$. 
Let $\alpha_i \in E^*$ such that $(\alpha_i,\theta)_{w_j} = 1$ if $j = 1,\dots,u_i$,
and that $(\alpha_i,\theta)_{w_j} = - 1$ if $j =  u_i+ 1, \dots , d-\tau_iv$.
\bigskip
Let $\ell : K \to k$ be a non--zero linear form such that $\ell(\overline x) = \ell(x)$ for
all $x \in K$; if ${\rm char}(K) = 0$, then we choose $\ell$ to be the trace map, $\ell = {\rm Tr}_{K/k} : K \to k$. 
Let $h' : V \times V \to K$ be the hermitian form defined by $h' = <\alpha_1,\dots,\alpha_m>$ and let $q_{h'}: V \times V \to k$ be
the quadratic space defined by $q_{h'}(x,y) = \ell (h'(x,y))$ for all $x,y \in V$.
By construction, the signature at $v$ of $q_{h'} $ is $(r_v,s_v)$. 

\bigskip
Let $S$ be the set of finite places of $k$ at which the Hasse invariants of $q$ and $q_{h'}$ are
not equal. This is a finite set of even cardinality : indeed, the Hasse invariants
of two quadratic spaces over $k$ differ at an even number of places, and 
$q$ and $q_{h'}$ are isomorphic at all the infinite places. 

\bigskip
Let $T$ be the set of finite places of $k$ such that every place of $E$ above $v$ splits
in $K$. 
Note that both quadratic spaces $q$ and $q_{h'}$ have isometries with minimal
polynomial $f$ over every completion of $k$ (by hypothesis for $q$, by construction
for $q_{h'}$). Therefore if $v \in T$, then both $q$ and $q_{h'}$ are hyperbolic 
over $k_v$ (cf. lemma 9.5). Hence $q$ and $q_{h'}$ are isomorphic over $k_v$,
and this implies that $v$ does not belong to $S$. Therefore  $S$ and $T$
are disjoint.
\bigskip
For each $v \in S$, let us choose a place $w$ of $E$ which does not split in $K$; this
is possible because $S$ and $T$ are disjoint. Let us denote by $S_E$ the set of these
places. Then $S_E$ is a finite set of even cardinality. 
\bigskip
For all $w \in S_E$, let $\beta_w \in E^*_w$ such that $(\beta_w,\theta)_w = -1$;
note that such a $\beta_w$ exists as
$w$ does not split in $K$.  By Hilbert's reciprocity, there exists $\beta \in E^*$
such that $(\beta,\theta)_w = (\beta_w,\theta)_w = -1$ if $w \in S_E$, and that
$(\beta_w,\theta)_w = 1$ for all the other places $w$ of $E$ (see for instance [O'M 73], 71:19, or
[PR 10], lemma 6.5). 
Let $h : V \times V \to K$ be the hermitian form given by 
$h = <\beta \alpha_1,\dots,\alpha_d>$ and let 
 $q_h : V \times V 
\to k$ be the quadratic space defined by $q_h (x,y) = \ell (h(x,y))$ for all
$x,y \in V$. 
Then 
by th. 7.2, the Hasse invariants of $q_h$ and  $q$ are equal. This implies that
 $q$ and  $q_h$ have equal dimension, determinant, signatures and Hasse
 invariants, therefore these quadratic spaces
are isomorphic. Note that $q_h$ has an isometry with minimal polynomial $f$ by
construction, hence $q$ also has such an isometry, and this concludes the proof
of the theorem.

\bigskip
{\bf \S 10. A necessary and sufficient condition}

\bigskip
Suppose that $k$ is a global field, and let us denote by $\Sigma_k$ the set of
all places of $k$. Let $q$ be a quadratic space over $k$,  and let $M$ 
be a module.  The aim of this section is to give some necessary and sufficient
conditions for $q$ to have an isometry with module $M$ (see th. 10.11 (b)).
This was already started in the previous section. One of the results of \S 9 can
be reformulated as follows :

\bigskip
\noindent
{\bf Theorem 10.1.} {\it Let $f \in k[X]$ be a symmetric, irreducible polynomial
of even degree {\rm (in other words, an irreducible polynomial of type 1)}. Let
$m \in {\bf N}$, let 
$M = [k[X]/(f)]^m$, and let $q$ be a quadratic space over $k$. Then $q$ has
an isometry with module $M$ if and only if such an isometry exists
over all the completions of $k$.}
\medskip
\noindent
{\bf Proof.} As $f$ is irreducible, a quadratic space $q$ of dimension $m {\rm deg}(f)$
has an isometry with minimal polynomial $f$ if and only if $q$ has an isometry
with module $M$. Hence the result follows from th. 9.1.
\bigskip

We have 
$M = M^0 \oplus M^1 \oplus M^2$ with $M_i$ of type $i$. Recall that a 
quadratic space has an isometry with a module of type 2 if and
only if it is hyperbolic (cf. prop. 3.5). Hence $q$ has an isometry with module $M$ if and
only it is isomorphic to an orthogonal sum of a quadratic space having
an isometry with module $M^0 \oplus M^1$ and of a hyperbolic space.
Therefore it suffices to consider modules with $M^2 = 0$.
\bigskip
On the other hand, we have seen that a quadratic space has an isometry with module
$M$ if and only if it is the orthogonal sum of a quadratic space having an isometry
with module $\overline M_{\rm odd}$ and of a hyperbolic space.  Since $\overline M_{\rm odd}$
is semi--simple, it is sufficient to treat the case of semi--simple modules.
\bigskip
Suppose that  $M$ is semi--simple, and that $M^2 = 0$. We have $M = M^0 \oplus M^1$, with $M^0 = [k[X]/(X+1)]^{n_+} \oplus [k[X]/(X-1)]^{n_-}$ for some $n_+, n_- \in {\bf N}$, and
$$M^1 \simeq [k[X]/(f_1)]^{n_1} \oplus \dots
\oplus  [k[X]/(f_r)]^{n_r},$$ where $f_1,\dots,f_r \in k[X]$ are distinct irreducible
polynomials of type 1 and $n_i \in {\bf N}$. Recall from \S 1 that this implies
that ${\rm deg}(f_i)$ is even, and that $f_i(1)f_i(-1) \not = 0$ for all $i = 1,\dots, r$.
Set $M_i =  [k[X]/(f_i)]^{n_i}$.
\bigskip
Set $M_0 = M^0$. Then we have $M=M_0 \oplus M_1 \oplus \dots \oplus M_r$. Set
 $I = \{ 1,\dots,r \}$, and $I_0 = I \cup \{0 \} =  \{ 0,\dots,r \}$.

\bigskip
\noindent
{\bf Proposition 10.2.} {\it  The quadratic space $q$ has an isometry
with module $M$ over $k$ if and only if there exist quadratic spaces
$q_0,\dots,q_r$ defined over $k$ such that $$q \simeq q_0 \oplus \dots \oplus q_r,$$
and that for all $i \in I_0$, the quadratic space $q_i$ has an isometry with module $M_i$
over all the completions of $k$.}
\bigskip
\noindent
{\bf Proof.} Suppose that $q$ has an isometry with module $M$ over $k$. Then by 3.3 there exist quadratic spaces $q_0,\dots,q_r$ defined over $k$ such that 
$q \simeq q_0 \oplus \dots \oplus q_r$
and that the quadratic space $q_i$ has an isometry with 
module $M_i$ over $k$ for all $i \in I_0$. 
\medskip
Let us prove the converse. By hypothesis  there exist quadratic spaces
$q_0,\dots,q_r$ defined over $k$ such that $q \simeq q_0 \oplus \dots \oplus q_r,$
and that for all $i \in I_0$, the quadratic space $q_i$ has an isometry with module $M_i$
over all the completions of $k$.
By th. 10.1 this
implies that the quadratic space $q_i$ has an isometry with module $M_i$
over $k$ for all $i \in I$. As $M_0$ is of type 0 and ${\rm dim}(q_0) = {\rm dim}(M_0)$, the
quadratic space $q_0$ has an isometry with module $M_0$. Since  $q$ is the orthogonal sum of the $q_i$'s, this proves the proposition. 

\bigskip
Suppose that $q$ has an isometry with module $M$ over all the completions of $k$.
Then by prop. 3.3 there exist quadratic spaces $\tilde q_i^v$ having an isometry with module $M_i$ for all $i \in I_0$ and for all
places $v$ of $k$ such that we have an isomorphism over $k_v$ 

$$q \simeq \tilde q^v_0 \oplus \dots \oplus \tilde q^v_r.$$

The quadratic spaces $\tilde q_i^v$  are not uniquely determined by $q$ and $M$. The strategy
used in this section is to investigate under what condition one can modify them
and obtain quadratic spaces $q^v_0,\dots,q^v_r$ defined over $k_v$ 
such that  that there exist quadratic spaces $q_0,\dots,q_r$ defined over $k$ 
which are isomorphic to $q^v_0,\dots, q^v_r$  over $k_v$ for all places
$v$ of $k$. By prop. 10.2, the Hasse principle holds precisely when this is possible. 
We start with some  definitions and  lemmas.
\bigskip
Let  us consider {\it collections} $C = \{ q^v_i \} $, for $i \in I_0$
and $v \in \Sigma_k$,  of quadratic spaces defined over $k_v$, and let us denote by  ${\cal C}_{M}$ the set of $C = \{ q^v_i \} $ of collections satisfying
the  condition 

\medskip
\noindent
(i) For all $v \in \Sigma_k$ and all $i \in I_0$, the quadratic space $q^v_i$ has an isometry with module $M_i$ over $k_v$.

\bigskip
Further, let us denote by  ${\cal C}_{M,q}$ the set of $C = \{ q^v_i \} \in {\cal C_M}$ of collections 
satisfying the additional  condition

\medskip
\noindent
(ii) For all  $v \in \Sigma_k$, we have $q \simeq  q^v_0 \oplus \dots \oplus q^v_r $
over $k_v$.

\bigskip
The above considerations show that if $q$ has an isometry with module $M$ over
$k_v$ for all the places $v$ of $k$, then there exist quadratic spaces $ q_i^v$ satisfying
(i) and (ii), in other words such that $C = \{ q^v_i \} \in {\cal C}_{M,q}$. Hence we have

\bigskip
\noindent
{\bf Lemma 10.3.} {\it Suppose that the quadratic space $q$ has an isometry
with module $M$ over $k_v$ for all  $v \in \Sigma_k$. Then ${\cal C}_{M,q}$ is
not empty.}

\bigskip
Note that if $C = \{ q^v_i \} \in {\cal C}_{M,q}$, then we have ${\rm dim}(q^v_i) = 
{\rm dim}(M_i)$ and ${\rm det}(q^v_i) = [f_i(1)f_i(-1)]^{n_i} \in k^*_v/k_v^{*2}$
for all places $v$ of $k$, and for all $i \in I$. Therefore the collections in 
${\cal C}_{M,q}$ can only differ by the Hasse invariants and the
signatures of the quadratic spaces. 

\bigskip
For any $C = \{ q^v_i \} \in {\cal C}$ and any $v \in \Sigma_k$, set

$$S_v(C) =  \{ i  \in I \ |  \ w(q^v_i) = 1 \} .$$ 

\bigskip

For all $i \in I$, set $d_i = [f_i(1)f_i(-1)]^{n_i}$, and let $d_0 = {\rm det}(q)d_1\dots d_r$. Set $D =  \Sigma_{i < j} (d_i,d_j) \in {\rm Br}_2(k)$. 
Recall that $w(q) \in  {\rm Br}_2(k)$ is the Hasse  invariant of $q$.  For any 
$x \in {\rm Br}_2(k)$ and any $v \in \Sigma_k$, let us denote by $x_v  \in \{0,1 \}$ the image of $x$ in ${\rm Br}_2(k_v)$.
\bigskip
\noindent
{\bf Proposition 10.4.} {\it Let $C = \{ q^v_i \} \in {\cal C}_{M}$. Then $C \in {\cal C}_{M,q}$ if and only if ${\rm det}(q^v_0) = d_0$, 
$$|S_v(C)| \equiv  w(q)_v  + D_v  \ \ ({\rm mod} \ 2)$$ for all $v \in \Sigma_k$,  and
${\rm sign}(q) = {\rm sign}(q_0^v \oplus \dots \oplus q_r^v) $ for all real places $v$ of $k$.}
\bigskip
\noindent
{\bf Proof.} Suppose that $C \in {\cal C}_{M,q}$. Then $q \simeq  q^v_0 \oplus \dots \oplus q^v_r $
for all $v \in \Sigma_k$. In particular, if $v$ is a real place, then 
${\rm sign}(q) = {\rm sign}(q_0^v \oplus \dots \oplus q_r^v) $. For all $v \in \Sigma_k$ and
all $i \in I$, 
we have ${\rm det}(q_i^v) = d_i$, hence ${\rm det}(q_0^v) = {\rm det}(q) d_1 \dots d_r = d_0$. Moreover, for all $v \in \Sigma_k$ we have 
$$w(q)_v = w(q_0^v \oplus \dots \oplus q_r^v) = w(q_0^v) + \dots + w(q_r^v) + \Sigma_{i < j} (d_i,d_j)
= |S_v(C)| + D_v,$$ as claimed.

\bigskip Let us prove the converse. 
Let $v \in \Sigma_k$ be a finite place, and let us check that  $q \simeq  q^v_0 \oplus \dots \oplus q^v_r $
over $k_v$. 
As  $C \in {\cal C}_{M}$,
the quadratic space $q^v_i$ has an isometry with module
$M_i$ over $k_v$ for all $i \in I_0$. Therefore
we have ${\rm det}(q^v_i) = d_i$ for all $i \in I$. By hypothesis, ${\rm det}(q^v_0) = d_0 =
{\rm det}(q) d_1\dots d_r$.
Thus
$$w(q_0^v \oplus \dots \oplus q_r^v) = w(q_0^v) + \dots + w(q_r^v) + \Sigma_{i < j} (d_i,d_j)
= |S_v(C)| + D  = w(q)_v.$$  Therefore $q_v$ and $q_0^v \oplus \dots \oplus q_r^v$ have
equal dimension, determinant and Hasse--Witt invariant, hence these quadratic spaces
are isomorphic over $k_v$.  If $v$ is a real place, then we are assuming that
${\rm sign}(q) = {\rm sign}(q_0^v \oplus \dots \oplus q_r^v) $, hence 
$q \simeq  q^v_0 \oplus \dots \oplus q^v_r $ over $k_v$. Thus
condition (ii) holds.  Since $C \in  {\cal C}_{M}$, condition (i) holds as well,
therefore we have $C \in {\cal C}_{M,q}$.

\bigskip
\noindent
{\bf Corollary 10.5.} {\it Let $\tilde C = \{ \tilde q^v_i \} \in {\cal C}_{M,q}$, and let $C = \{ q^v_i \} \in {\cal C}_{M}$.
Let $u \in \Sigma_k$ be a finite place, and let $\alpha, \beta \in I_0$
with $\alpha \not = \beta$ such that 
\medskip
{\rm (a)} $q^v_i \simeq \tilde q^v_i$ for all $v \not = u$ and for all $i \in I_0$;

\medskip
{\rm (b)} $q^u_i \simeq \tilde q^u_i$ for all $i \not = \alpha, \beta$;

\medskip
{\rm (c)} $w(q^u_{\alpha})  \not = w(\tilde q^u_{\alpha}) $ and $w(q^u_{\beta})  \not = w(\tilde q^u_{\beta}) $;

\medskip
{\rm (d)} ${\rm det}(q_0) = {\rm det}(q) d_1 \dots d_r$.

\medskip
Then $C \in {\cal C}_{M,q}$.}

\bigskip
\noindent
{\bf Proof.}  By (b) and (c), we have $|S_u(C)| = |S_u(\tilde C)|$. By {\rm (a)}, we have
$|S_v(C) = |S_v(\tilde C)|$ for all  $v \not = u$, and 
${\rm sign}(q_0^v \oplus \dots \oplus q_r^v)  = {\rm sign}(\tilde q^v_0 \oplus \dots \oplus \tilde q^v_r)$,
if $v$ is a real place. Since  $\tilde C \in {\cal C}_{M,q}$, by prop. 10.4 we have $|S_v(\tilde C)| \equiv w(q)_v + D_v \ \ ({\rm mod} \ 2)$  for all $v \in \Sigma_k$, and ${\rm sign}(q) = {\rm sign}(\tilde q_1^v \oplus \dots \oplus \tilde q_r^v) $
if $v$ is a real place. Hence we also have $|S_v(C)| \equiv w(q)_v  + D_v \ \ ({\rm mod} \ 2)$
for all $v \in \Sigma_k$, and ${\rm sign}(q) = {\rm sign}(q^v_0 \oplus \dots \oplus q^v_r) $
for all real places $v$ of $k$. By prop. 10.4, this implies that $C  \in {\cal C}_{M,q}$. 

\bigskip
\noindent
{\bf Lemma 10.6.} {\it Let $v$ be a finite, non--dyadic place of $k$, let $i \in I_0$, and let $Q$ a quadratic space over $k_v$
with module $M_i$. }
\medskip
\noindent
 {\rm (a)} {\it Suppose that $i \not = 0$. Then there exists
 a quadratic space $Q'$ over $k_v$ having an isometry with module $M_i$ such
 that $w(Q') = 0$.}
 \medskip
 \noindent
 {\rm (b)} {\it Suppose that $i = 0$, and let $d \in k^*/k^{*2}$. Then there exists a
 quadratic space $Q'$ over $k_v$ having an isometry with module $M_0$ such
 that $w(Q') = 0$ and ${\rm det}(Q') = d$.}
 
 \bigskip
 \noindent
 {\bf Proof.}  {\rm (a)} If $w(Q) = 0$, there is nothing to prove. Suppose that
 $w(Q) = 1$. Since $v$ is non--dyadic, this implies that the quadratic space $Q$ is
 not hyperbolic. Therefore by cor. 3.5 the module $M_i$ is not hyperbolic over $k_v$. By
lemma 9.6, there exists a quadratic space $Q'$ over $k_v$ having an isometry with
module $M_i$ such that $w(Q') = 0$.
\medskip
\noindent
{\rm (b)} Set $n_0 = {\rm dim}(M_0)$, and let $Q'$ be the $n_0$--dimensional
quadratic space $$Q' = <1,\dots,1,d>.$$ Then ${\rm det}(Q') = d$ and $w(Q') = 0$.
As any quadratic space of dimension $n_0$ has
an isometry with module $M_0$, this completes the proof of the lemma.

\bigskip In order to give a necessary and sufficient condition for the Hasse principle
to hold, 
the first step is to show that ${\cal C}_{M,q}$ contains 
a collection $C = \{ q^v_i \} $ in ${\cal C}_{M,q}$ such that $w(q^v_i) = 0$
for almost all places $v$ of $k$ and all $i \in I_0$. Recall that
$D =  \Sigma_{i < j} (d_i,d_j) \in {\rm Br}_2(k)$. Let $T$ be the set of places $v$ of $k$ such that $D_v \not = 0$, and let
$S$ be the set of places of $k$ at which the Hasse invariant of $q$ is
not equal to the Hasse invariant of the hyperbolic space of dimension equal to ${\rm dim}(q)$. 
Let $\Sigma_2$ be the set of dyadic places and $\Sigma_{\infty}$ the set of infinite places
of $k$. 
Set $\Sigma = S \cup T \cup \Sigma_2 \cup \Sigma_{\infty}$. Note that $\Sigma$ is
a finite subset of $\Sigma_k$.

\bigskip
\noindent
{\bf  Proposition 10.7.} {\it The set ${\cal C}_{M,q}$ contains a collection $C = \{ q^v_i \} $ of quadratic forms defined over $k_v$ such that $w(q^v_i) = 0$
for all $v \not \in \Sigma$  and all $i \in I_0$.}

\bigskip
\noindent
{\bf Proof.} Let $\tilde C = \{ \tilde q^v_i  \} \in {\cal C}_{M,q}$. 
Let $v$ be a place of $k$ such that $v \not \in \Sigma$ and suppose that
$|S_v(\tilde C)| \not = 0$. It suffices to show that there exists a collection $C \in {\cal C}_{M,q}$
with $|S_v(C)| <  |S_v(\tilde C)|$.
\medskip
 Set $w^v_i = w(\tilde q^v_i)$.  We are supposing that $|S_v(\tilde C)| \not = 0$, hence 
there exists an $i$ with  $w^v_i = 1$.
Since $v \not \in S \cup \Sigma_2$, we have $w(q)_v = 0$. Moreover $v \not \in T$, hence
$w(q)_v  = w^v_0  + \dots + w^v_r $. Thus there exists $j \not = i$ such
that $w^v_j = 1$. 
\medskip
By lemma 10.6 there
exist quadratic spaces 
$q_i^v$  and $q_j^v $ over $k_v$ having isometries
with module $M_i$ respectively $M_j$ such that ${\rm det}(q_i^v) = d_i$, ${\rm det}(q_j^v) = d_j$
and 
$w(q^v_i) = w(q^v_j) = 0$. Set $q^v_{\alpha} =  \tilde q^v_{\alpha}$ if $\alpha \not = i,j$. 

\medskip
Set  $C = \{  q^v_i \}$. Then $C$ satisfies the conditions of cor. 10.5, hence
 $C \in {\cal C}_{M,q}$. Note that
$C = \{  q^v_i \}$ satisfies $w(q^v_i) = w(q^v_j) = 0$, therefore $|S_v(C)| <  |S_v(\tilde C)|$.
This completes the proof of the proposition.

\bigskip
 For any collection $C =  \{ q^v_i \} \in {\cal C}_{M,q}$ and all $i \in I_0$, set
 $$T_i (C) =  \{ v  \in \Sigma_k \ |  \ w(q^v_i) = 1 \} .$$
Let ${\cal F}_{M,q}$ be the subset of ${\cal C}_{M,q}$ consisting of the collections
$C = \{ q^v_i \} $ of quadratic spaces over $k_v$ such that for all $i \in I_0$, the set
$T_i(C)$ is finite. 

\bigskip
\noindent
{\bf Theorem 10.8.} {\it Suppose that $q$ has an isometry with module $M$ over
$k_v$ for all places $v$ of $k$. Then $q$ has an isometry with module $M$ if 
and only if there exists
a collection $C = \{ q^v_i \} \in {\cal F}_{M,q}$ such that for all $i \in I_0$, the cardinality of
$T_i(C)$ is even.}
\bigskip
\noindent
{\bf Proof.} Suppose that $q$ has an isometry with module $M$. Then by prop. 10.2, 
there exist quadratic spaces
$q_0,\dots,q_r$ defined over $k$ such that $q \simeq q_0 \oplus \dots \oplus q_r,$
and that the quadratic space $q_i$ has an isometry with module $M_i$
over all the completions of $k$  for all $i \in I_0$. Let $q_i^v = q_i \otimes_k  k_v$, and
let $C = \{ q^v_i \}$. Then $C \in {\cal C}_{q,M}$, and for all $i = 0,\dots,r$,  the set
$T_i(C)$ is finite of even cardinality.
\bigskip
Conversely, let $C = \{ q^v_i \} \in {\cal F}_{M,q}$ be such that 
$T_i(C)$ has
even cardinal for all $i \in I_0$. Recall that as the quadratic space $q^v_i$ has an isometry with module $M_i$,
we have ${\rm dim}(q^v_i) = 
{\rm dim}(M_i)$ and ${\rm det}(q^v_i) = d_i \in k^*_v/k_v^{*2}$
for all places $v$ of $k$, and for all $i \in I_0$. Therefore by [O'M 73], Chapter VII, th. 72.1,
for all $i \in I_0$ there exists a quadratic space $q_i$ such that $q_i \otimes_k k_v
\simeq q_i^v$ for all $v \in \Sigma_k$. We have $q \simeq q_0 \oplus \dots \oplus q_r$
over $k_v$ for all $v$, hence by the Hasse--Minkowski theorem we have
$q \simeq q_0 \oplus \dots \oplus q_r.$ Therefore by prop. 10.1, the quadratic
space $q$ has an isometry with module $M$.

 \bigskip
 For any module $N$ and any $d \in k^*$,  let $\Omega(N,d)$ be the
 set of finite places $v$ of $k$  such that for any $\epsilon \in \{0,1 \}$, there exists a
 quadratic space $Q$ over $k_v$ with ${\rm disc}(Q) = d$ and $w(Q) = \epsilon$
 having an isometry with module $N \otimes_k k_v$.

 \bigskip
 For all $i,j \in I_0$, let $\Omega_{i,j} = \Omega (M_i,d_i) \cap \Omega (M_j,d_j)$.
 \bigskip
 \noindent
 {\bf Remark 10.9.} Note that if $i,j \in I$, then $\Omega_{i,j}$ does not depend on $q$. If
 $M_0 \not = 0$ and $i = 0$, then $\Omega_{i,j}$ depends on $d_0 = {\rm det}(q) d_1 \dots d_r$.
 \bigskip
  Recall that 
 for any collection $C =  \{ q^v_i \} \in {\cal F}_{M,q}$ and all $i \in I_0$, we have
 $$T_i (C) =  \{ v  \in \Sigma_k \ |  \ w(q^v_i) = 1 \}.$$
 \medskip
 \noindent
 {\bf Definition 10.10.} We say that $C = (q_i^v) \in {\cal F}_{M,q}$ is {\it connected} if 
for all $i \in I$ such that  $|T_i(C)|$ is odd,
there exist $j \in I$ with $j \not = i$ such that  $|T_j(C)|$  is odd,
and  a chain $i = i_1, \dots, i_m = j$ of elements of $I$ with  $\Omega_{i_t,i_{t+1}}  \not = \emptyset$ for all $t = 1, \dots, m-1$.
 We say that  ${\cal F}_{M,q}$ is
{\it connected} if it contains a connected element. 

\bigskip
\noindent
{\bf Theorem 10.11.} (a) {\it The quadratic space $q$ has an isometry with module $M$ over
$k_v$ for all $v \in \Sigma_k$ if and only if ${\cal F}_{M,q}$ is not empty.}
\medskip
\noindent
(b) {\it The quadratic space $q$ has an isometry with module $M$ over
$k$ if and only if ${\cal F}_{M,q}$ is connected.}

\medskip
\noindent
{\bf Proof.} {\rm (a)}  It is clear that if ${\cal F}_{M,q}$ not empty, the quadratic space $q$ has an isometry
with module $M$ over $k_v$ all $v \in \Sigma_k$. The 
converse follows from lemma 10.3 and th. 10.7.
\medskip \noindent
{\rm (b)}  If the quadratic space $q$ has an isometry with module $M$, then there
exist quadratic spaces $q_0,\dots,q_r$ over $k$ such that $q \simeq q_0 \oplus \dots \oplus q_r$
and that $q_i$ has an isometry with module $M_i$ for all $i \in I_0$. 
Set $q^v_i = q_i \otimes_k k_v$, and let $C = (q^v_i)$. Then $C  \in {\cal F}_{M,q}$,
and $|T_i(C)|$ is even for all $i \in I_0$. Therefore $C$ is a connected element of ${\cal F}_{M,q}$, hence ${\cal F}_{M,q}$ is connected.
\medskip
Conversely, suppose that ${\cal F}_{M,q}$ is connected, and let $C =  (q_i^v) \in {\cal F}_{M,q}$ be a connected element.
Suppose that for some $i \in I_0$, the integer $|T_i(C)|$ is odd. Since $C$ is connected, there exist
$j \in I_0$ with $j \not = i$ such that  $|T_j(C)|$  is odd,
and  a chain $i = i_1, \dots, i_m = j$ of elements of $I$ with  $\Omega_{i_t,i_{t+1}} \not = \emptyset$ for all $t = 1, \dots, m-1$.  Let $v_t  \in \Omega_{i_t,i_{t+1}}$.
Then there exist  quadratic spaces $\tilde q^{v_t}_t $ over $k_v$ with  
$w(\tilde q^{v_t}_t)  \not = w(q^{v_{t}}_{t})$ and ${\rm det}(\tilde q^{v_t}_t) = d_t$ having an isometry with module $M_t$. Set 
$\tilde q^u_s = q^u_s$ if $(u,s) \not = (v_t,t)$.  Set $\tilde C = (\tilde q_i^v)$; then
$\tilde C  \in {\cal F}_{M,q}$. We have $|T_i (\tilde C)|
 \equiv  0  \ \ ({\rm mod} \ 2)$, $|T_j (\tilde C)|
 \equiv  0  \ \ ({\rm mod} \ 2)$, and $|T_s (\tilde C)|
 \equiv  |T_s (C)|  \ \ ({\rm mod} \ 2)$ if $s \not = i,j$. 
Repeating this procedure we obtain
a family of quadratic spaces $C'  \in {\cal F}_{M,q}$ such that $|T_i(C)|$ is even for all $i \in I_0$.
By th. 10.8 this implies that $q$ has an isometry with module $M$.

\medskip
Note that the condition (a) does not imply condition (b) in general  (in other words, there
are counter-examples to the Hasse principle) :  this follows from
the examples of Prasad and Rapinchuk, cf. [PR 10], Example 7.5.

\bigskip
 \bigskip
 {\bf \S 11. The case of modules of mixed type}
 
 \bigskip
 The aim of this section and the next is to give some applications of th. 10.11. We keep the notation
 of the previous section : in particular, $k$ is a global field and $\Sigma_k$ is the set
 of places of $k$. Recall that $M$ is semi-simple, and that $M \simeq M^0 \oplus M^1$,
 with $M^0$ of type 0 and $M^1$ of type 1. If $M^1 = 0$, then we already have a
 complete criterion for the existence of an isometry with module $M$ (see prop. 4.3).
 In this section, we consider the case where both $M^0$ and $M^1$ are non--zero.
 As we will see, the case where 
 ${\rm dim}(M^0) \ge 3$ is especially simple, and will be considered first. Then
 we examine the case where ${\rm dim}(M^0) = 2$ or $1$.
Let $q$ be
 a quadratic space over $k$, and assume that ${\rm dim}(q) = {\rm dim}(M)$.
 
 \bigskip
 \noindent
 {\bf Definition 11.1.} 
 For every real place $v$ of $k$, let $(r_v,s_v)$ denote the signature of $q$ over $k_v$,
and let $\sigma_v$ be the number of roots of $F_M \in k_v[X]$ that are not on the unit circle. 
We say that the {\it signature conditions} are  satisfied if for every real place 
$v$ of $k$, we have $(r_v,s_v) \ge (\sigma_v,\sigma_v)$, and if moreover $M^0 = 0$, then 
$(r_v,s_v) \equiv (\sigma_v,\sigma_v)
\ ({\rm mod} \ 2)$. 

\bigskip
\noindent
{\bf Proposition 11.2.} {\it Suppose that ${\rm dim}(M^0) \ge 3$.
 Then the quadratic space $q$ has an isometry with module $M$ if and
 only if the signature conditions are satisfied.}
 
 \bigskip
 The proof of prop. 11.2, as well that of several other results of sections 11 and 12, is based
 on prop. 11.3 below.  With the notation of \S 10, we have :

 \bigskip
 \noindent
 {\bf Proposition 11.3.} {\it Suppose that there exists $i_0 \in I_0$ such that for all
 $i \in I_0$ we have $\Omega_{i_0,i}(q)  \not = \emptyset$. Suppose that the quadratic space $q$
 has an isometry with module $M$ over every completion of $k$. Then $q$ has an
 isometry with module $M$. }

\bigskip
For the proof of prop. 11.3, we need the following lemmas. We use the notation of \S 10.
\bigskip
\noindent
{\bf Lemma 11.4.} {\it Let $C \in {\cal F}_{M,q}$. Then  $$\Sigma_{v \in \Sigma_k} | S_v(C)| 
 \equiv  0  \ \ ({\rm mod} \ 2).$$}
 
 \bigskip
 \noindent
 {\bf Proof.} By prop. 10.4, we have
 $$|S_v(C)| \equiv  w(q)_v  + D_v  \ \ ({\rm mod} \ 2)$$ for all $v \in \Sigma_k$.
Hence 

$$\Sigma_{v \in \Sigma_k} | S_v(C)| 
 \equiv  \Sigma_{v \in \Sigma_k}  w(q)_v +  \Sigma_{v \in \Sigma_k} D_v \ \ ({\rm mod} \ 2).$$
 As $w(q)$ and $D$ are elements of ${\rm Br}_2(k)$, we have  $$\Sigma_{v \in \Sigma_k}  w(q)_v 
  \equiv  0  \ \ ({\rm mod} \ 2), \ \ {\rm and} \ \ 
\Sigma_{v \in \Sigma_k}  D_v   \equiv  0  \ \ ({\rm mod} \ 2).$$ This implies that
$\Sigma_{v \in \Sigma_k} | S_v(C)| 
 \equiv  0  \ \ ({\rm mod} \ 2)$, as claimed.

 \bigskip
 \noindent
 {\bf Lemma 11.5.}  {\it Let $C \in {\cal F}_{M,q}$. Then  $$\Sigma_{i \in I_0}  | T_i (C)| 
 \equiv  0  \ \ ({\rm mod} \ 2).$$}

 \noindent
 {\bf Proof.} Indeed, we have $$\Sigma_{i \in I_0}  | T_i (C)|  = \Sigma_{v \in \Sigma_k} | S_v(C)| .$$
 By lemma 11.4, we have $\Sigma_{v \in \Sigma_k} | S_v(C)| 
 \equiv  0  \ \ ({\rm mod} \ 2)$, hence $\Sigma_{i \in I}  | T_i (C)| 
 \equiv  0  \ \ ({\rm mod} \ 2).$

\bigskip
\noindent
{\bf Proof of Proposition 11.3.} By th. 10.11 (a), the set ${\cal F}_{M,q}$ is not empty. Let $C =  (q_i^v) \in {\cal F}_{M,q}$, and let  $i \in I_0$
be such that $|T_i(C)|$ is odd. By lemma 11.5, we have $\Sigma_{i \in I_0}  | T_i (C)| 
 \equiv  0  \ \ ({\rm mod} \ 2)$, hence there exists $j \in I_0$ with $j \not = i$ such
 that $|T_i(C)|$ is odd. By hypothesis, we have $\Omega_{i,i_0} \not = \emptyset$
 and $\Omega_{j,i_0} \not = \emptyset$, hence $C$ is connected. Therefore ${\cal F}_{M,q}$
 is connected, hence by th. 10.11 (b) the quadratic space $q$ has an isometry with
 module $M$.

 \bigskip
 \noindent
 {\bf Lemma 11.6.} {\it Let $N$ be a module of type 0, and let $d \in k^*/k^{*2}$. 
 \medskip
 \noindent
{\rm (a)} If ${\rm dim}(N) \ge 3$, then every  finite place of $k$ is in $\Omega (N,d)$.
\medskip
\noindent
{\rm (b)} If ${\rm dim}(N) = 2$ and $d \not = -1 \in k^*/k^{*2}$, then every finite place of $k$ is
in $\Omega (N,d)$.}

\medskip
\noindent
{\bf Proof.} Since $N$ is of type 0, every quadratic space of dimension 
equal to ${\rm dim}(N)$
 has an isometry with module $N$. Therefore the result follows from [O'M 73, 63:23].
 
 \bigskip
 \noindent
 {\bf Proposition 11.7.} {\it Suppose that ${\rm dim}(M^0) \ge 3$, or
 ${\rm dim}(M^0) = 2$ and ${\rm det}(q) \not = - d_1 \dots d_r \in k^*/k^{*2}$. 
 If the quadratic space $q$ has an isometry with module $M$ over every
 completion of $k$, then $q$ has an isometry with module $M$ over $k$.}
 \medskip
 \noindent
 {\bf Proof.}  If ${\rm dim}(M^0) \ge 3$, then lemma 11.6 (a) implies that every
 finite place of $k$ is in $\Omega(M^0,d_0)$. Therefore $\Omega_{0,i} \not = \emptyset$ 
 for all $i \in I_0$. Suppose that ${\rm dim}(M^0) = 2$ and ${\rm det}(q) \not = - d_1 \dots d_r \in k^*/k^{*2}$.  Recall that $d_0 = {\rm det}(q) d_1 \dots d_r \in k^*/k^{*2}$.  Hence 
 $d_0 \not = -1 \in k^*/k^{*2}$, and therefore by lemma 11.6 (b) every finite place of
 $k$ is in $\Omega(M^0,d_0)$. This implies that $\Omega_{0,i} \not = \emptyset$ 
 for all $i \in I_0$ in this case as well, and hence 
 the proposition  follows from prop. 11.3.

 \bigskip
 \noindent
 {\bf Proof of prop. 11.2.} The necessity of the signature conditions  follows from th. 8.1. Let us
 show that they are also sufficient. By prop. 11.7, it suffices to show that $q$ has an isometry with module $M$
 over $k_v$ for all $v \in \Sigma_k$. For real places, this is a consequence of th. 8.1. 
 Let $v$ be a finite place. If $M^1 = 0$, then $M$ is of type 0, and every quadratic
 space of dimension equal to ${\rm dim}(M)$ has an isometry of module $M$. Suppose
 that $M^1 \not = 0$, and note that this implies that ${\rm dim}(M^1) \ge 2$.  We have $M^1 \otimes _k k_v \simeq N^v_1 \oplus N^v_2$
 where $N^v_1$ is of type 1 and $N^v_2$ of type 2. If $N^v_1 \not = 0$, then
 the result follows from th. 7.1. Suppose that $N^v_1 = 0$, and let $2m_2 = {\rm dim}(N^v_2)$.
Since ${\rm dim}(M^1) \ge 2$, we have ${\rm dim}(M) \ge 5$, hence $q$ is isotropic over $k_v$,
and its  Witt index is $\ge m_2$. By prop. 4.3, this implies that $q$ has an
isometry with module $M$ over $k_v$. This concludes the proof of the
proposition. 
\medskip
\noindent
{\bf Proposition 11.9.} {\it Suppose that ${\rm dim}(M^0) = 2$ and that
${\rm det}(q) \not  = - d_1 \dots d_r \in k^*/k^{*2}$. Then $q$ has an isometry
with module $M$ if and only  the following two conditions hold :}
\medskip
\noindent
(a){\it  The signature conditions are satisfied}; 
\medskip
\noindent
(b)  {\it If $v$ is a finite place and if $M^1 \otimes_k k_v$ is hyperbolic, then the
Witt index of $q$ over $k_v$  is 
$\ge {1 \over 2}  {\rm dim}(M^1)$. }

\medskip
\noindent
{\bf Proof.} By prop. 11.7, we have to show that the conditions hold if and only if
$q$ has an isometry with module $M$ over $k_v$ for all $v$. For real places,
this is a consequence of th. 8.1. Let $v$ be a finite place. Recall  that $M^1 \not = 0$, and let
$M^1 \otimes _k k_v \simeq N^v_1 \oplus N^v_2$
 where $N^v_1$ is of type 1 and $N^v_2$ of type 2. If $N^v_1 \not = 0$, then
 the result follows from th. 7.1. Suppose that $N^v_1 = 0$, and 
note that  this means that $M^1 \otimes _k k_v$ is hyperbolic. 
By prop. 4.3, this implies that $q$ has an
isometry with module $M$ over $k_v$ if and only if the Witt index of $q$ over $k_v$ is
$\ge {1 \over 2} {\rm dim}(M^1)$, and this is precisely condition (b). 
\bigskip
\noindent
{\bf Proposition 11.10.} {\it Suppose that ${\rm dim}(M^0) = 2$ and that
${\rm det}(q) = - d_1 \dots d_r \in k^*/k^{*2}$. Then $q$ has an isometry
with module $M$ if and only if $q \simeq q_0 \oplus q'$ where $q_0$ is
a hyperbolic plane, and $q'$ is a quadratic space over $k$ having
an isometry with module $M^1$.}
\medskip
\noindent
{\bf Proof.} If $q \simeq q_0 \oplus q'$ with $q_0$ a hyperbolic plane and
$q'$ a quadratic space having an isometry with module $M^1$, then
$q$ has an isometry with module $M$. Let us prove the converse.
Suppose that $q$ has an isometry with module $M$. Then $q \simeq q_0 \oplus q'$
with $q_0$ a quadratic space having an isometry with module $M_0$, and $q'$
a quadratic space having an isometry with module $M^1$. As $M^1$ is of
type 1, we have ${\rm det}(q') = d_1\dots d_r \in k^*/k^{*2}$. By hypothesis,
we have ${\rm det}(q) = - d_1 \dots d_r \in k^*/k^{*2}$. Therefore
${\rm det}(q_0) = -1 \in k^*/k^{*2}$. Since ${\rm dim}(q_0) = 2$, this implies
that $q_0$ is isomorphic to a hyperbolic plane. 
\bigskip
Recall that $d_0 = {\rm det}(q) d_1 \dots d_r \in k^*/k^{*2}$.
\bigskip
\noindent
{\bf Proposition 11.11.} {\it Suppose that ${\rm dim}(M^0) = 1$. Then $q$ has
an isometry with module $M$ if and only if $q \simeq q_0 \oplus q'$ where
$q_0 \simeq <d_0>$ and $q'$ is a quadratic space having an isometry
with module $M^1$.}
 \medskip
 \noindent
 {\bf Proof.}  If $q \simeq q_0 \oplus q'$ with $q_0 \simeq <d_0>$  and
$q'$ a quadratic space having an isometry with module $M^1$, then
 $q$ has an isometry with module $M$. Let us prove the converse.
Suppose that $q$ has an isometry with module $M$. Then $q \simeq q_0 \oplus q'$
with $q_0$ a quadratic space having an isometry with module $M_0$, and $q'$
a quadratic space having an isometry with module $M^1$. As $M^1$ is of
type 1, we have ${\rm det}(q') = d_1\dots d_r \in k^*/k^{*2}$. Since
${\rm dim}(q_0) = {\rm dim}(M_0) = 1$ and $d_0 = {\rm det}(q) d_1 \dots d_r \in k^*/k^{*2}$,
we have  $q_0 \simeq <d_0>$. This concludes the proof of the proposition.

\bigskip
\bigskip
{\bf \S 12. Modules of type 1}

\bigskip
We keep the notation of sections 10 and 11. In particular, $k$ is a global
field and $\Sigma_k$ is the set of places of $k$. In this section, we assume
that $M$ is a semi-simple module of type 1. Recall that this means that 
$M \simeq M_1 \oplus \dots \oplus M_r$, where 
$M_i = [k[X]/(f_i)]^{n_i}$ for some symmetric,
irreducible polynomials  $f _i \in k[X]$ of even degree, and for some $n_i \in {\bf N}$.
We use the notation  $I = \{1, \dots, r \}$,
and $K_i = k[X]/(f_i)$. Let $q$ be a quadratic space over $k$ such that
${\rm dim}(q) = {\rm dim}(M)$. Recall that we denote by $F_M \in k[X]$ the characteristic polynomial of $M$. If $v$ is a real place of $k$, then we denote
by $(r_v,s_v)$ the signature of $q$ at $v$, and by $\sigma_v$ the number
of roots of $F_M \in k^*/k^{*2}$ off the unit circle. 
\bigskip
Recall that the {\it signature conditions} are satisfied for $q$ and $M$ if for all real
places $v$ of $k$, we have $(r_v,s_v) \ge (\sigma_v,\sigma_v)$, and  
$(r_v,s_v) \equiv (\sigma_v,\sigma_v)
\ ({\rm mod} \ 2)$. 
\bigskip
We say that the {\it hyperbolicity conditions} are satisfied for $q$ and $M$ if for all places
$v$ of $k$ such that $M \otimes_kk_v$ is a hyperbolic module (that is, a module of type 2),  the quadratic form
$q_v$ over $k_v$ is hyperbolic. 

\bigskip 
We have the following
\bigskip
\noindent
{\bf Theorem 12.1.} {\it The quadratic space $q$ has an isometry with module $M$ over
all the completions of $k$ if and only if the signature conditions and the hyperbolicity
conditions are satisfied, and we have ${\rm det}(q) = F_M(1)F_M(-1) \in k^*/k^{*2}.$}
\medskip
{\bf Proof.} This follows from cor. 3.5, cor. 7.4, and  th. 8.1.

\bigskip
We will see that the necessary and sufficient conditions of th. 10.11 can be
interpreted in terms of splitting properties of the fields $K_i$. We start with a few lemmas.
Let us recall that 
 for any module $N$ and any $d \in k^*$,  we denote by  $\Omega(N,d)$ the
 set of finite places $v$ of $k$  such that for any $\epsilon \in \{0,1 \}$, there exists a
 quadratic space $Q$ over $k_v$ with ${\rm disc}(Q) = d$ and $w(Q) = \epsilon$
 having an isometry with module $N$. 
 \bigskip
 \noindent
 {\bf Lemma 12.2.} {\it Let $f \in k[X]$ be a symmetric, irreducible polynomial of even degree, and let
 $m \in {\bf N}$.  Set 
 $N = [k[X]/(f)]^m$ and let $d = [f(1)f(-1)]^m$. Let $v$ be a finite place of $k$. Then $v \in \Omega (N,d)$
 if and only if $N\otimes_k k_v$ is not hyperbolic.}
 \medskip
 \noindent
 {\bf Proof.} If $N\otimes_k k_v$ is hyperbolic, then every quadratic space with module
 $N\otimes_k k_v$ is hyperbolic (cf. lemma 3.5.), therefore $v \not  \in \Omega (N,d)$. 
 Conversely, suppose that $N\otimes_k k_v$ is not hyperbolic. Then by lemma 9.6,
 for any $\epsilon \in \{0,1 \}$ there exists a quadratic space $Q$ having an isometry
 with module $N \otimes_k k_v$ such that $w(Q) = \epsilon$. By cor. 5.2, we
 have ${\rm det}(Q) = d$, hence $v \in  \Omega (N,d)$.

 \bigskip
 \noindent
 {\bf Notation 12.3. } Let $E$ be an extension of finite degree of $k$, let $K$ be a quadratic
 extension of $E$, and let  $x \mapsto \overline x$ be the non--trivial
 automorphism of $K$ over $E$. Let us denote by $\Sigma^s(K)$ the set of
 $v \in \Sigma_k$ such that every place of $E$ above $v$ splits in $K$. Let
 $\Sigma^{ns}(K)$ be complement of $\Sigma^s(K)$ in $\Sigma_k$; in other words,
 $v \in \Sigma_k$ such that there exists a place of $E$ above $v$ that is not split
 in $K$.

 \bigskip
 Let $\Sigma_k'$ be the set of finite places of $k$. Then we have
 \bigskip
 \noindent
 {\bf Lemma 12.4.} {\it For all $i \in I$, we have $\Omega(M_i,d_i) = \Sigma^{ns}(K_i) \cap \Sigma'_k$.}
 
 \medskip
 \noindent
 {\bf Proof.} Let $v \in \Sigma'_k$. By lemma 9.4, we have $v \in \Sigma^{ns}(K_i)$ if and only if $M_i \otimes_k k_v$ is
 not hyperbolic, and by lemma 12.2 this is equivalent to $v \in \Omega (M_i,d_i)$. 
 
 \bigskip
 For all $i,j \in I$, set $\Sigma^{ns}_{i,j} = \Sigma^{ns}(K_i) \cap \Sigma^{ns}(K_j)$. 
 
 \bigskip
 \noindent
 {\bf Theorem 12.5.} {\it Assume that there exists $i_0 \in I$ such that
 for all $i \in I$, we have $\Sigma^{ns}_{i_0,i} \not = \emptyset$. Suppose that $q$ has
 an isometry with module $M$ over all the completions of $k$. Then $q$ has an
 isometry with module $M$. }

 \medskip
 \noindent
 {\bf Proof.} Let $i \in I$, and let us show that there exists a finite place $v$ of $k$ such that $v \in \Sigma^{ns}_{i_0,i}$. 
Indeed, let $u$ be a real place of $k$ with $u  \in \Sigma^{ns}_{i_0,i}$. 
Let $L$ be a Galois extension of $k$ containing the fields $K_{i_0}$ and $K_i$, and let
$G = {\rm Gal}(L/k)$. Let us denote by $c$ the conjugacy
class of the complex conjugation in $G$ corresponding to an extension of
the place $u$ to $L$. By the Chebotarev density theorem,
there exists a finite place $v$ of $k$ such that the conjugacy class
of the Frobenius automorphism at $v$ is equal to $c$. Let $v$ be such a place.
Then all the places of $E_{i_0}$, respectively  $E_i$, above $v$ are inert in $K_{i_0}$,
respectively $K_i$. Therefore, 
we have $v \in \Sigma^{ns}_{i_0,i} =  \Sigma^{ns}(K_{i_0}) \cap
\Sigma^{ns}(K_i)$. 
Since $v \in \Sigma'_k$, by lemma 12.4  this implies that $v \in 
\Omega (M_{i_0},d_{i_0}) \cap \Omega (M_i,d_i) = \Omega_{i_0,i}$. Therefore, for all
$i \in I$, we have $\Omega_{i_0,i} \not = \emptyset$. Since here $I = I_0$, prop. 11.3 gives
the desired result. 

\bigskip
\noindent
{\bf Corollary 12.6.}  {\it Suppose  that there exists $i_0 \in I$ such that
 for all $i \in I$, we have $\Sigma^{ns}_{i_0,i} \not = \emptyset$.   Then $q$ has an
isometry with module $M$ if and only if the hyperbolicity and signature conditions
are satisfied, and if we have ${\rm det}(q) = F_M(1)F_M(-1) \in k^*/k^{*2}.$}
\medskip
\noindent
{\bf Proof.} This follows from th. 12.1 and prop. 12.5.
\bigskip
The hypothesis of th. 12.5 and cor. 12.6  are often satisfied, for instance we have

\bigskip
\noindent
{\bf Corollary 12.7.} {\it Suppose that there exists a real place $v$ of $k$ such that
all the roots of $F_M \in k_v[X]$ are on the unit circle. Then $q$ has an
isometry with module $M$ if and only if the hyperbolicity and signature conditions
are satisfied, and if we have ${\rm det}(q) = F_M(1)F_M(-1) \in k^*/k^{*2}.$}

\medskip
\noindent
{\bf Proof.} Indeed, we have $v \in \Sigma^{ns}(K_i)$ for all $i \in I$, hence $v \in \Sigma^{ns}_{i,j}$ for all $i,j$. Therefore the result follows from cor. 12.6. 

\bigskip
Recall that a number field is CM if it is a totally imaginary quadratic extension
of a totally real number field. We say that the module $M$ is of {\it type CM} if
$k = {\bf Q}$ and if  the fields $K_i$ are CM fields for all $i \in I$. 

\bigskip
\noindent
{\bf Corollary 12.8.} {\it Suppose that $M$ is a module of type CM. Then the
quadratic space has an isometry with module $M$ if and only if the hyperbolicity
conditions are satisfied, if ${\rm det}(q) = F_M(1)F_M(-1) \in k^*/k^{*2}$, and
the signature of $q$ is even.}

\medskip
\noindent
{\bf Proof.} Indeed, the hypothesis of cor. 12.7 is satisfied, and the signature
condition amounts to saying that the signature of $q$ is even. 
\bigskip
\bigskip
{\bf Bibliography}
\bigskip 
\noindent
[B 84] E. Bayer--Fluckiger, Definite unimodular lattices having an automorphism of given characteristic polynomial.  {\it Comment. Math. Helv.}  {\bf 59}  (1984),  509--538.
\medskip
\noindent
[B 99] E. Bayer--Fluckiger, Lattices and number fields.  Algebraic geometry: Hirzebruch 70 (Warsaw, 1998),  69--84, {\it Contemp. Math.}  {\bf 241}, Amer. Math. Soc., Providence, RI, 1999. 
\medskip 
\noindent
[B 12] E. Bayer--Fluckiger, Embeddings of tori in orthogonal groups, {\it Ann. Inst. Fourier},
to appear.
\medskip
\noindent
[BMa 13] E. Bayer-Fluckiger, P. Maciak, Euclidean minima of abelian fields of odd prime power conductor,  {\it Math. Ann.}, to appear.
\medskip 
\noindent
[BMar 94] E. Bayer-Fluckiger and J. Martinet, R\'eseaux
li\'es aux alg\`ebres
semi-simples, {\it J. reine angew. Math.} {\bf 451}
(1994), 51--69.
\medskip 
\noindent
[BCM 03] R. Brusamarello, P. Chuard--Koulmann and J. Morales, Orthogonal groups containing
a given maximal torus, {\it J. Algebra} {\bf 266} (2003), 87--101. 
\medskip 
\noindent
[F 12] A. Fiori,  Special points on orthogonal symmetric spaces, {\it J. Algebra}, {\bf 372} (2012),
397-419.
\medskip 
\noindent
[GR 12]  S. Garibaldi and A. Rapinchuk, Weakly commensurable S--arithmetic subgroups
in almost simple algebraic groups of types B and C, {\it Algebra and Number Theory},
to appear. 
\medskip 
\noindent
[G 05] P. Gille, Type des tores maximaux des groupes semi-simples, {\it J. Ramanujan
Math. Soc.}  {\bf 19} (2004), 213--230. 
\medskip
\noindent
[GM 02] B. Gross and C.T. McMullen, Automorphisms of even unimodular lattices and unramified Salem numbers.  {\it J. Algebra} {\bf  257 }  (2002),  265--290. 
\medskip
\noindent
[JRV 12] F. Jouve and F. Rodriguez--Villegas, On the bilinear structure associated to the Bezoutians, preprint 2012.
\medskip 
\noindent
[Lee 12] T-Y. Lee, Embedding functors and their arithmetic properties, {\it Comment. Math. Helv.}, to appear.
\medskip
\noindent
[Le 69] J. Levine, Invariants of knot cobordism,  {\it  Invent. Math} {\bf 8}  (1969), 98--110.
\medskip
\noindent
[Le 80] J. Levine, {\it Algebraic structure of knot modules}, Lecture Notes in Mathematics {\bf 772}, Springer Verlag, Berlin, 1980.
\medskip 
\noindent
[M 69] J. Milnor, Isometries of inner product spaces, {\it  Invent. Math.} {\bf  8}  (1969), 83--97. 
\medskip 
\noindent
[O'M 73] O.T. O'Meara, {\it Introduction to quadratic forms}, Reprint of the 1973 edition. Classics in Mathematics. Springer-Verlag, Berlin, 2000.
\medskip
\noindent
[PR 10] G. Prasad and A.S. Rapinchuk, Local--global principles for embedding of fields with
involution into simple algebras with involution, {\it Comment. Math. Helv.} {\bf 85} (2010),
583--645. 
\medskip
 \noindent
[Sch 85] W. Scharlau, {\it Quadratic and hermitian forms}, Grundlehren der Mathematischen Wissenschaften {\bf 270}, Springer-Verlag, Berlin, 1985.

\bigskip
\bigskip
Eva Bayer--Fluckiger

EPFL-FSB-MATHGEOM-CSAG

Station 8

1015 Lausanne, Switzerland

eva.bayer@epfl.ch

\bye